\documentclass[10pt]{article}

\usepackage{a4wide}
\usepackage{amssymb}
\usepackage{amsfonts}
\usepackage{amsmath}
\input xy
\xyoption{arrow} \xyoption{matrix}

\date{}

\newtheorem{proposition}{Proposition}[section]
\newtheorem{theorem}[proposition]{Theorem}
\newtheorem{lemma}[proposition]{Lemma}

\newtheorem{corollary}[proposition]{Corollary}

\def\Kdim{{\rm K.dim }\,}

\def\der{\partial }

\def\nFM0{{\nu }_{F,M_0}}
\def\nFN0{{\nu }_{F,N_0}}
\def\nGN0{{\nu }_{G,N_0}}

\def\N0{ {\bf N}_0 }

\def\t{\otimes}
\def\g{\gamma}
\def\v{\varphi}
\def\ra{\rightarrow}

\def\Xpm{X^{\pm }}

\def\s{\sigma}
\def\Z{{\bf Z }}

\def\l1{{\lambda}_1}

\def\a{\alpha}
\def\a0{ {\alpha }_0}
\def\a1{ {\alpha }_1}

\def\l{\lambda}
\def\o{\omega}

\def\nFGM0{{\nu }_{F,G,M_0}}

\def\nFN0{{\nu}_{F,N_0}}

\def\sm{{\sigma}^m}

\def\sm1{{\sigma}^{-1}}

\def\smtp1{{\sigma}^{-t+1}}

\def\o{\omega }
\def\S1{S^{-1}}

\def\Xpm1{X^{\pm 1}_1}

\def\sPM1{{\sigma }^{\pm 1}}
\def\sMP1{{\sigma }^{\mp 1 }}

\def\d{\delta}

\def\di{{\rm d.ind}}

\def\L{\Lambda}

\def\CA{{\cal A}}

\def\CD{{\cal D}}

\def\Ytm1{Y^{t-1}}
\def\Yim1{Y^{i-1}}

\def\k{{\bf k}}

\def\Aut{{\rm Aut}}
\def\bK{\overline{K}}

\def\Der{{\rm Der }}
\def\ad{{\rm ad }}
\def\dim{{\rm dim }}

\def\gcd{ {\rm gcd } }

\def\SL2Z{ {\rm SL}_2({\bf Z}) }

\def\CR{ {\cal R}}

\def\Gp1{ G^{1 , 1 } }
\def\P11{ P^{-1 , 1 } }
\def\Pp1{ P^{1 , 1 } }

\def\Supp{{\rm Supp}}

\def\nCLsr{{}^\nu\kern-2pt {\cal L}^{\sigma , \rho  }}
\def\nP{{}^\nu \kern-2pt P}
\def\nL{{}^\nu\kern-2pt L}
\def\nLL{{}^\nu\kern-2pt \Lambda}
\def\nPsr{{}^\nu\kern-2pt P^{\sigma , \rho  }}
\def\nLsr{{}^\nu\kern-2pt L^{\sigma , \rho  }}
\def\nuCL{{}^\nu\kern-2pt  {\cal L}}
\def\nCLsr{{}^\nu\kern-2pt {\cal L}^{\sigma , \rho  }}
\def\nCL1m{{}^\nu\kern-2pt {\cal L}^{-1 , 1  }}
\def\x1nu{x^\frac{1}{\nu}}
\def\xm1nu{x^{-\frac{1}{\nu}}}

\def\CR{ {\cal R}}

\def\ra{\rightarrow }

\def\CB{{\cal B}}

\def\CC{ {\cal C}}

\def\nAM0{{\nu }_{{\cal A},M_0}}
\def\nAN0{{\nu }_{{\cal A},N_0}}

\def\Kdim{ {\rm Kdim } }

\def\End{ {\rm End }}
\def\Der{ {\rm Der }}

\def\CR{ {\cal R }}

\def\ad{ {\rm ad }}

\def\ga{\mathfrak{a}}

\def\gm{\mathfrak{m}}
\def\gp{\mathfrak{p}}
\def\gq{\mathfrak{q}}

\def\gM{\mathfrak{M}}

\def\Spec{{\rm Spec}}

\def\di!{\frac{\der^i}{i!}}
\def\dik!{\frac{\der^k_i}{k!}}

\def\CC{{\cal C}}

\def\Max{{\rm Max}}

\def\Fp{\mathbb{F}_p}

\def\Max{{\rm Max}}

\def\gldim{{\rm gldim}}

\def\N{\mathbb{N}}

\def\0{\overline{0}}
\def\1{\overline{1}}

\def\Ln1{\L_{n,\overline{1}}}

\def\a1{a_{\overline{1}}}

\def\St{{\rm St}}
\def\S{\Sigma}

\def\vn1{\overrightarrow{n-1}}

\def\Sh{{\rm Sh}}

\def\im{{\rm im}}

\def\F{\mathbb{F}}

\def\mS{\mathbb{S}}
\def\mJ{\mathbb{J}}
\def\mI{\mathbb{I}}

\def\ann{{\rm ann}}

\def\mF{\mathbb{F}}

\def\mT{\mathbb{T}}

\def\mE{\mathbb{E}}

\def\mU{\mathbb{U}}

\def\K1{{\rm K}_1}

\def\Supp{{\rm Supp}}

\def\hmI1{\widehat{\mI_1}}
\def\tmI1{\widetilde{\mI_1}}
\def\tmJ1{\widetilde{\mJ_1}}
\def\hB1{\widehat{B_1}}
\def\hCB1{\widehat{\CB_1}}

\def\ga{\mathfrak{a}}

\def\oG{\overline{G}}

\def\Z{\mathbb{Z}}

\def\CB{{\cal B}}
\def\CC{{\cal C}}

\def\hB{\hat{B}}

\def\FFp{\overline{\F}_p}
\def\Fr{{\rm Fr}}

\def\Prim{{\rm Prim}}

\begin{document}

\author{V. V.   Bavula 
} 
 
\title{Isomorphism problems and groups  of automorphisms for  Ore extensions $K[x][y; f\frac{d}{dx} ]$ (prime characteristic)}

\maketitle 
\begin{abstract} 
Let $\L (f) = K[x][y; f\frac{d}{dx} ]$ be an Ore extension of a  polynomial algebra $K[x]$ over an arbitrary field $K$ of characteristic $p>0$  where $f\in K[x]$.  For each polynomial $f$,  the automorphism group of the algebras $\L (f)$ is explicitly described. The automorphism group $\Aut_K(\L (f))=\mS\rtimes G_f$ is a semidirect product  of two explicit groups where $G_f$ is the {\em eigengroup} of the polynomial $f$ (the set of all automorphisms of $K[x]$ such that $f$ is their common eigenvector).  For each polynomial $f$, the eigengroup $G_f$ is explicitly described. It is proven that every subgroup of $\Aut_K(K[x])$ is the eigengroup of a polynomial.

It is proven that the Krull and global dimensions of the algebra $\L (f)$ are 2. The prime, completely prime, primitive and maximal ideals of the algebra $\L (f)$ are classified. \\



 {\em Mathematics subject classification
2010:  16D60, 13N10, 16S32, 16P90,  16U20.

$${\bf Contents}$$
\begin{enumerate}
\item Introduction.
\item Spectra, the centre, the Krull and global dimensions of the algebra $\L $. 
\item Isomorphism problems and groups  of automorphisms for  Ore extensions $K[x][y; f\frac{d}{dx} ]$. 
\item The eigengroup $G_f$ of a polynomial $f$. 

\end{enumerate}
}

\end{abstract}


\section{Introduction}

In this paper, module means a left module, $K$ is a field of characteristic $p>0$ and $\bK$ is  its algebraic closure, $K^\times := K\backslash \{ 0\}$,  $K[x]$ be a polynomial algebra in the variable  $x$ over $K$, $\Der_K(K[x])=K[x]\frac{d}{dx}$ is the set of all $K$-derivations of the algebra $K[x]$, 
$$\L :=\L (f):=K[x][y; \d := f\frac{d}{dx} ]=K\langle x, y \, |  \, yx-xy=f\rangle =\bigoplus_{i\geq 0} K[x]y^i$$
is an Ore extension of the algebra $K[x]$ where $f=f(x)\in K[x]$. Given  an algebra $D$ and its derivation $\d$, the {\em Ore extension} of $D$, denoted $D[y; \d]$,  is an algebra generated by the algebra $D$ and $y$ subject to the defining relations $yd-dy=\d (d)$ for all $d\in D$.  The algebra $\L$ is a Noetherian domain of Gelfand-Kirillov dimension 2.\\

The aim of the paper is for each polynomial $f$ to give an explicit description of the  automorphism group $\Aut_K(\L (f))$ of the algebra $\L (f)$.\\

By dividing the element $y$ by the leading coefficient of the polynomial $f$ we can assume that the polynomial $f$ is {\em monic}, i.e. its leading coefficient is 1 provided $f\neq 0$. Then the  algebras $\{ \L (f)\, | \, f\in K[x]\}$ as a class can be partitioned into four subclasses: $f=0$, $f=1$, the polynomial $f$ has only a {\em single} root in $\bK$ and the polynomial $f$ has at least two {\em distinct} roots in $\bK$.\\

 If $f=0$ then the algebra $\L (0)=K[x,y]$ is a polynomial algebra in two variables and its group of automorphisms is well-known \cite{van der Kulk-1953}:
  The group $\Aut_K (K[x,y])$ is generated by the automorphisms:
 \begin{eqnarray*}
 t_l &:&  x\mapsto \l x, \;\;  \;\; \;\; \;\;  \;\; \; y\mapsto y, \\
\Phi_{n,\l} &:& x\mapsto x+\l y^n, \;\; y\mapsto y,  \\
 \Phi_{n,\l}' &:& x\mapsto x, \;\;  \;\; \;\; \;\;  \;\; \;\; \;\; y\mapsto y+\l x^n,
\end{eqnarray*}
 where $n\geq 0$ and $\l \in K$. \\

  If $f=1$ 
   then 
 the algebra $\L (1)$ is the (first) {\em Weyl algebra} $$A_1=K\langle 
x, \der\, | \, \der x-x\der =1\rangle\simeq K[x][y; \frac{d}{dx} ].$$
 In characteristic zero  Dixmier  \cite{Dix}, and in prime characteristic   Makar-Limanov \cite{Makar-Limanov-1984}, gave an explicit generators for  the automorphism group $\Aut_K (A_1)$ (see also \cite{Bav-AutWeylCharp} for more results on $\Aut_K (A_1)$): The group $\Aut_K (A_1)$ is generated by the automorphisms:
 \begin{eqnarray*}
\Phi_{n,\l} &:& x\mapsto x+\l y^n, \;\; y\mapsto y,  \\
 \Phi_{n,\l}' &:& x\mapsto x, \;\, \; \;\; \;\;  \;\; \;\; \;\; y\mapsto y+\l x^n,
\end{eqnarray*}
 where $n\geq 0 $ and $\l \in K$. \\
 
 Extending results of Dixmier \cite{Dix}  on the automorphisms of the Weyl algebra $A_1$,  
Bavula and Jordan \cite{Bav-Jor} considered isomorphisms and automorphisms of generalized
Weyl algebras over polynomial algebras of characteristic 0. Explicit generators of the automorphism group were given (they are more involved comparing to the case of the Weyl algebra to present them here). Alev and Dumas \cite{Alev-Dumas-1997}
initiated the study of automorphisms of Ore extensions $\L (f)$ in characteristic zero case. Their results were extended also to prime characteristic by Benkart, Lopes and Ondrus \cite{Benkart-Lopes-Ondrus-2015}. The  algebra $\L (x^2)$ (the Jordan plane)  was studied by Shirikov \cite{Shirikov-2015},  Cibils, Lauve, and  \cite{CLW}, and  Iyudu \cite{Iyudu-2014}. Gadis \cite{Gaddis-2015} studied isomorphism problems for  algebras on two generators that satisfy a single quadratic relation.\\
 
 {\bf Isomorphism problems  for  the algebras $\L (f)$.} Theorems \ref{Alev-Dum-P3.6}  is an isomorphism criterion for the algebras $\L$. It also describes the automorphism group of each algebra $\L (f)$.

 \begin{theorem}\label{Alev-Dum-P3.6}
 Let $f, g\in K[x]$ be   polynomials. Then $\L (f) \simeq \L(g)$ iff $g(x)= \l f(\alpha x+\beta)$ for some elements $\l , \alpha \in K^\times$ and $\beta \in K$.
\end{theorem}
  
In characteristic zero,  Theorem \ref{Alev-Dum-P3.6} was proven by Alev and Dumas \cite[Proposition 3.6]{Alev-Dumas-1997}  (1997) and in prime characteristic -- by Benkart, Lopes and Ondrus 
\cite[Theorem 8.2]{Benkart-Lopes-Ondrus-2015} (2015).

Benkart, Lopes and Ondrus  \cite{Benkart-Lopes-Ondrus-2015} (Theorems 8.3 and 8.6) gave a 
description of the {\em set} of  automorphisms groups  of algebras $\L (f)$ over arbitrary fields and if the  automorphism  group of $\L (f)$ is {\em given} they presented  information  on the type of  the polynomial $f$, \cite[Corollary 8.7]{Benkart-Lopes-Ondrus-2015} (in general, if one fixes the type of the polynomial then the automorphism group is {\em larger} than the one which  is naively expected). In this paper, we proceed in the opposite direction: if the polynomial $f$ is {\em given} then the automorphism group $\Aut_K\, \L (f)$   is explicitly described. \\

{\bf The eigengroup $G_f(K)$ of a polynomial $f\in K[x]$.} Recall that $\Aut_K(K[x])=\{ \s_{\l , \mu} \, | \, \l \in K^\times, \mu \in K\}$  where $\s_{\l , \mu } (x) = \l x+\mu$.\\

{\em Definition, \cite{Bav-AutOreExt}.}  For a polynomial $f\in K[x]$, 
\begin{equation}\label{GfDef}
G_f=G_f(K):=\{ \s \in \Aut_K(K[x]) \, | \, \s (f)=\l_\s f\;\; {\rm for\; some}\;\; \l_\s\in K^\times\}
\end{equation}
 is called the {\em eigengroup} of the polynomial $f$. \\

Clearly, the set $G_f(K)$ is a subgroup of  $\Aut_K(K[x])$, it is the  largest subgroup of  $\Aut_K(K[x])$ such that the polynomial $f$  is their common eigenvector. For all scalars $\mu \in K$, $G_\mu =\Aut_K(K[x])$. For  all scalars  $\nu\in K^\times$, $G_f=G_{\nu f}$. So, in order to describe the eigengroup $G_f(K)$ we can assume that the polynomial $f$ is a  monic  polynomial.
 It is proven that every subgroup of $\Aut_K(K[x])$ is the eigengroup of a polynomial (Theorem \ref{30Mar20}). For each subgroup $G$ of $\Aut_K(K[x])$ all the polynomials $f\in K[x]$ with $G_f=G$ are explicitly described in the case when the field $K$ is algebraically closed. The most interesting and difficult case is when the group $G$ is a finite group. There are three types of finite subgroups in $\Aut_K(K[x])$ that are not the identity group. For each such group $G$, the polynomials $f$ with $G_f=G$ has a unique form/presentation the, so-called. {\em eigenform} of $f$.

 The eigengroup $G_f$ has an isomorphic copy in the automorphism group $\Aut_K(\L (f))$: The map 
\begin{equation}\label{GfAutL}
G_f(K)\ra \Aut_K(\L (f)), \s_{\l  \mu} \mapsto \s_{\l , \mu}: x\mapsto \l x+\mu , \;\; y\mapsto y\mapsto \l^{\deg (f) -1} y
\end{equation}
is a group monomorphism where $\deg (f)$ is the degree of the polynomial $f$. We identify the group $G_f(K)$ with its image in $\Aut_K(\L (f))$. The group $G_f(K)$ is the most important/difficult part of the group $\Aut_K(\L (f))$. Approximately half of the paper is about how to find it. If the group $G_f(K)$ is a finite group  it is a semidirect product of two subgroups, $G_f(K)=\widetilde{G}_f(K)\rtimes\oG_f(K)$ (Theorem \ref{C7Mar20}).\\

There are four distinguish cases:
\begin{enumerate}
\item $\widetilde{G_f}\neq \{ e\}$, $\oG_f\neq \{ e\}$,
\item $\widetilde{G_f}\neq \{ e\}$, $\oG_f=\{ e\}$,
\item $\widetilde{G_f}=\{ e\}$, $\oG_f\neq \{ e\}$,
\item $\widetilde{G_f}= \{ e\}$, $\oG_f=\{ e\}$.
\end{enumerate}
In the case when  $K=\bK$, Theorem \ref{A10Mar20}, Theorem \ref{A11Mar20}, Theorem \ref{B11Mar20} and Theorem \ref{C11Mar20} are  criteria  for each case to hold, respectively (see also Proposition \ref{A9Mar20}). These four theorems 
 are also  explicit descriptions of the eigengroup $G_f(K)$. They also show that in each of four cases the polynomial $f$ admits a {\em unique} presentation -- {\em the eigenform} of $f$ (introduced in the paper).

 At the end of Section \ref{pAUTR}, a finite  algorithm is given of finding the eigengroups $G_f(\bK)$ and $G_f(K)$,  and the eigenform of $f$. \\
 
 In the case when $K\neq \bK$, similar results are obtained, see  Theorem \ref{A12Mar20}. Proposition \ref{a26Mar20} gives criteria for the groups $\widetilde{G}_f(K)$, $\oG_f(K)$ and $G_f(K)$ to be $\{ e\}$. \\
 
 {\bf The group of automorphisms of  the algebra $\L (f)$.}  
 Given a group $G$, a normal subgroup $N$ and a subgroup $H$. The  group $G$ is called the {\em semidirect product} of $N$ and $H$, written $G=N\rtimes H$,  if $G=NH:=\{ nh\, | \, n\in N, h\in H\}$ and $N\cap H=\{ e\}$ where $e$ is the identity of the group  $G$. 
 
  The automorphism group $\Aut_K(\L (f))$ of the algebra $\L (f)$ contains an obvious subgroup
\begin{equation}\label{grmS}
\mS :=\mS (K):= \{ s_p\, | \, p\in K[x]\}\simeq (K[x],+), \;  s_p\mapsto p\;\; {\rm where}\;\;  s_p(x) =x\; {\rm and } \;  s_p(y)=y+p.
\end{equation}

\begin{theorem}\label{p5Mar20}
 Suppose that $f\in K[x]$ is a monic nonscalar polynomial. Then $$\Aut_K(\L (f))=\mS (K) \rtimes G_f(K).$$ 
\end{theorem}



{\bf The Krull and global dimensions of the algebra $\L (f)$.} It it proven that the  Krull and global dimensions of the algebra $\L (f)$ are 2 (Theorem \ref{A16Mar20}).\\

{\bf Classifications of prime, completely prime, primitive and maximal ideals  of the algebra $\L (f)$.} Theorem \ref{pKGLD} classifies  prime, completely prime, primitive and maximal ideals  of the algebra $\L (f)$.  The height 1 prime ideals were classified in \cite[Theorem 7.6]{Benkart-Lopes-Ondrus-2015}. It is proven that every nonzero ideal of the algebra $\L (f)$ meets the centre of $\L (f)$ (Corollary \ref{a16Mar20}).  \\

{\bf Classifications of simple $\L (f)$-modules.} In \cite{Bav-SimOre-1999}, simple modules were classified for all Ore extensions $A=D[x; \s , \der ]$ where $D$ is 
a Dedekind domain, $\s\in $ is an automorphism of $D$ and $\der$ is a $\s$-derivation of $D$ (for all $a,b\in D$, $\der (ab)=\der (a)b+\s (a) \der (b)$).  Recall that the ring $A$ is generated by $D$ and $x$ subject to the defining relations: For all elements $d\in D$, $xd=\s (d)x+\der (d)$. 
 The algebras $\L (f)$ is a a very special case of the rings $A$.

Theorem \ref{17Mar20} classifies simple {\em left} $\L (f)$-modules, see also \cite{Benkart-Lopes-Ondrus-SimMod2015}.  For each simple left $\L (f)$-module an explicit $K$-basis is given and the actions of the canonical  generators $x$ and $y$ of the algebra $\L (f)$ on the basis is explicitly described. 

A classification of simple {\em right} $\L (f)$-modules is obtained at once from the classification of simple left $\L (f)$-modules by using the fact that the opposite algebra $\L(f)^{op}$ of the algebra $\L(f)$ is isomorphic to 
\begin{equation}\label{Lfop}
\L(f)^{op}\simeq \L(-f).
\end{equation}
Recall that the {\em opposite algebra} $A^{op}$ of an algebra $A$ coincides with the algebra $A$ as vector space but the multiplication in $A^{op}$  is given by the rule $a\cdot b=ba$. 
Every right $A$-module is a left $A^{op}$-module and vice versa.


\section{Spectra, the centre, the Krull and global dimensions of the algebra $\L $}\label{KGLDIM}

 In this section, $K$ is a field of characteristic $p>0$ (not necessarily algebraically closed) and $f=p_1^{n_1}\cdots p_s^{n_s}$ is a nonscalar polynomial of $ K[x]$ where $p_1, \ldots , p_s$ are irreducible, co-prime divisors of $f$ (i.e. $K[x]p_i+K[x]p_j=K[x]$ for all $i\neq j$).  The aim of this section is to find the centre of the algebra $\L (f)$;  to classify  simple $\L (f)$-modules;  to classify prime, completely prime, primitive and maximal ideals of the algebra $\L (f)$; and to prove that the Krull and global dimension of the algebra $\L (f)$ is 2. \\

{\bf The centre of the algebra $\L (f)$.} It follows from the direct sum
$$A_1=\bigoplus_{i,j=0}^{p-1}K[x^p, \der^p]x^i\der^j$$
and the commutation relation $[\der , x]=1$, that the centre $Z(A_1)$ of the Weyl algebra is equal to $K[x^p, \der^p]$, a polynomial algebra in two variables.  Let $K(x^p, \der^p)$ be the field of fractions of the polynomial algebra  $K[x^p, \der^p]$. 
The localization  $\CA_1$ of the Weyl algebra $A_1$ by the Ore set 
$K[x^p, \der^p]\backslash \{ 0\}$ is a simple $p^2$-dimensional algebra 
$$\CA_1=\bigoplus_{i,j=0}^{p-1}K(x^p, \der^p)x^i\der^j.$$
This follows from the relation $[\der , x]=1$. So, $Z(\CA_1)=K(x^p, \der^p)$. Hence, {\em every nonzero ideal of the Weyl algebra $A_1$ meets the centre of} $A_1$. 

The polynomial algebra $K[x]$ is a left $A_1$-module where $\der$ acts as the derivation $\frac{d}{dx}$. Furthermore, the kernel of the corresponding algebra homomorphism 
\begin{equation}\label{}
A_1\ra \End_K(K[x]), \;\; x\mapsto x, \;\; \der \mapsto \frac{d}{dx}
\end{equation}
is generated by the central element $\der^p$. So, $A_1/(\der^p)=\bigoplus_{i=0}^{p-1}K[x]\der^i\subseteq  \End_K(K[x])$.  The factor algebra $A_1/(\der^p)$ is a subalgebra of the algebra $\CD (K[x])$ of differential operators on the polynomial algebra $K[x]$ and the Weyl algebra $A_1$ is not. This is in sharp contrast with the characteristic zero case where $A_1=\CD (K[x])$. 

The algebra $\L =\L (f)$ can be identified with a subalgebra of  the Weyl algebra $A_1$  by the monomorphism:
\begin{equation}\label{(2.1)}
\L \ra A_1,\;\;x\mapsto x,\;\; y\mapsto f\der.
\end{equation}
So, $\L = K\langle  x, y=f\der\rangle\subset A_1$. Theorem \ref{16Mar20} describes the centre of the algebra $\L (f)$. It also gives explicit expressions for the $p$'th power of various elements of the algebras $\L (f)$ and $A_1$ that are key facts in finding the centre of $\L(f)$. 

The fact that the centre of the algebra $\L (f)$ is equal to $K[x^p, y^p-(\d^{p-2}(f))'y]$ was proven by Benkart, Lopes and Ondrus, \cite[Theorem 5.3,(2)]{Benkart-Lopes-Ondrus-2015}. Here we present a short proof of this fact. 

\begin{theorem}\label{16Mar20}
Let $\d = f\frac{d}{dx}\in \Der_K(K[x])$ where $f\in K[x]\backslash \{ 0\}$, and   $ g':=\frac{d g}{dx}$ where $g\in K[x]$. Then 
\begin{enumerate}
\item  $\d^p=(\d^{p-2}(f))'\d\in \Der_K(K[x])$. 
\item In the algebra $\L (f)$, $y^p=f^p\der^p+(\d^{p-2}(f))'y$. In particular, in the Weyl algebra $A_1$, $(f\der )^p=f^p\der^p+(\d^{p-2}(f))'f\der $.
\item The centre  $Z(\L (f))$ of the algebra $\L (f)$ is the polynomial algebra $$K[x^p, y^p-(\d^{p-2}(f))'y]=K[x^p,f^p\der^p]$$ and $f^p\der^p=y^p-(\d^{p-2}(f))'y$.
\item The algebra $\L (f) = \bigoplus_{i,j=0}^{p-1}Z(\L (f))x^iy^j$ is a free $Z(\L (f))$-module of rank $p^2$. 
\item The localization of the algebra $\L (f)$ at the Ore set $Z(\L (f))\backslash \{ 0\}$ is $\CA_1$.
\end{enumerate}
\end{theorem}

{\it Proof}. 1. Since $\d^p\in \Der_K(K[x])$, we have that $\d^p=g\frac{d}{dx}$ where  $g=\d^p (x)=\d^{p-1}(f)= (\d^{p-2}(f))'f$. Therefore, 
 $$\d^p= (\d^{p-2}(f))'f\frac{d}{dx}=(\d^{p-2}(f))'\d .$$ 
 
2. Notice that  $y^p=(f\der )^p=f^p\der^p+\sum_{i=1}^{p-1}a_i\der^i$ for some elements $a_i\in K[x]$. Recall that  $A_1/(\der^p)=\bigoplus_{i=0}^{p-1}K[x]\der^i\subseteq  \End_K(K[x])$. By statement 1, 
$$ (f\der )^p\equiv \sum_{i=1}^{p-1}a_i\der^i\equiv (\d^{p-2}(f))'\d \equiv (\d^{p-2}(f))'f\der \equiv (\d^{p-2}(f))'y \mod (\der^p),$$
and so  $a_1=(\d^{p-2}(f))'f$ and $a_2=\cdots = a_{p-1}=0$. 

3--5. By statement 2, $f^p\der^p=y^p-(\d^{p-2}(f))'y\in Z(\L (f))$. Hence, $$Z':=K[x^p, y^p-(\d^{p-2}(f))'y]=K[x^p,f^p\der^p]\subseteq Z(\L (f))$$ and $ \L (f) = \bigoplus_{i,j=0}^{p-1}Z'x^iy^j$. Now,
$$(Z'\backslash \{ 0\})^{-1}\L (f)= \bigoplus_{i,j=0}^{p-1}K(x^p, \der^p)x^iy^j=\bigoplus_{i,j=0}^{p-1}K(x^p, \der^p)x^i\der^j=\CA_1,$$
and so statement 5 is obvious and statements 3--4 follow. $\Box $

Let $A$ be an algebra and $a\in A$. The map $\ad_a =[a, -]: A\ra A$, $b\mapsto [a,b]:=ab-ba$ is a derivation of the algebra $A$ which is called the {\em inner derivation} of $A$ associated with the element $a$. 

\begin{corollary}\label{a16Mar20}
Every nonzero ideal of the algebra $\L (f)$ meets the centre of $\L (f)$.
\end{corollary}

{\it Proof}. Let $I$ be a nonzero ideal of the algebra $\L (f)$. Fix a nonzero element of $I$, say $a=\sum_{i,j=0}^{p-1}z_{ij}x^i\der^j$ for some elements $z_{ij}\in Z(\L (f))$,  by Theorem \ref{16Mar20}.(4). 
 Then applying several times the inner derivation $\ad_x :=[x, -]$ of the algebra $\L (f)$, we obtain a nonzero element, say $b\in I\cap Z(\L (f)[x]$. Then $0\neq b^p\in I\cap Z(\L (f))$. $\Box $\\

 {\bf The prime, completely prime, primitive and maximal  spectra of the algebra $\L $.}  An ideal $\gp$ of a ring $R$ is called a {\em completely prime ideal} if the factor ring $R/\gp$ is a domain. A completely prime ideal is a prime ideal. The sets of prime and completely prime ideals of the ring $R$ are denoted by $\Spec (R)$ and $\Spec_c(R)$, respectively.
The annihilator of a simple $R$-module is called a {\em primitive idea}l of $R$. Every primitive ideal is a prime ideal of $R$. The set of all primitive ideals is  denoted by $\Prim (R)$.  The set of all maximal ideals of $R$ is denoted by $\Max (R)$. Clearly, $\Max (R) \subseteq \Prim (R)\subseteq \Spec (R)$.

  An element $a$ of an algebra $A$ is called a {\em normal element} of $A$ if $Aa=aA$. An element $a$ of an algebra $A$ is called a {\em regular element} if it is not a zero divisor. The set of all regular elements of the algebra $A$ is denoted by $\CC_A$. Each regular normal element $a$ of the algebra $A$ determines an automorphism of the algebra $A$ given by the rule:
\begin{equation}\label{oaaut}
\o_a: A\ra A, \;\; b\mapsto \o_a(b)\;\; {\rm where}\;\; ab=\o_a(b)a.
\end{equation}
The elements $p_1, \ldots , p_s$ are regular normal elements of the algebra $\L =\L (f)$ (recall that $f=\prod_{i=1}^sp_i^{n_i}$)  since $$yp_i=p_i(y-p_i^{-1}f)\; \; {\rm and}\;\; xp_i=p_ix.$$
Therefore, $\o_{p_i} (x)=x$ and $\o_{p_i}(y) = y+p_i^{-1}f$.
 
 For an ideal $\ga$ of an algebra $A$, we denote by $V(\ga)$ the set of all prime ideals of $A$ that contain the ideal $\ga$. Let   $\min  \ga $  be  the {\em set of minimal primes} of $\ga$. These are the minimal elements of the set $V(\ga )$ with respect to inclusion. Suppose that the set $S_a:=\{ a^i\, | \, i\geq 0\}$ is a left Ore set of a domain $A$.  The algebra $A_a:=S_a^{-1}A=\{ a^{-i}b \, | \, b\in A\}$ is called the localization of $A$ at the powers of the element $a$. 
 
 For a commutative algebra $C$ and its non-nilpotent element $s$,
 the map 
 $$\Spec (C)\backslash V(s)\ra \Spec (C_s), \;\; \gp \mapsto S_s^{-1}\gp$$
  is a bijection with the inverse map $\gq\mapsto C\cap \gq$. We identify the sets $\Spec (C)\backslash V(s)$ and $ \Spec (C_s)$ via the bijection above, i.e.  $\Spec (C)\backslash V(s)= \Spec (C_s)$.
  
  By Theorem \ref{16Mar20}.(4), the centre of the Weyl algebra $A_1=K[x][\der ; \frac{d}{dx}]$ is the polynomial algebra $K[x^p, \der^p]$ (the result is well-known). Then $$f^p\in K[x^p]\subseteq Z(\L (f))=K[x^p, f^p\der^p]\subseteq K[x^p, \der^p]=Z(A_1)$$ and
\begin{equation}\label{Lffloc}
L=\L (f) \subseteq L_f=A_{1,f}= L_{f^p}=A_{1,f^p}=\bigoplus_{i,j=0}^{p-1} Z(L)_{f^p}x^iy^i=\bigoplus_{i,j=0}^{p-1}Z(A_1)_{f^p}x^i\der^i=\bigoplus_{i,j=0}^{p-1}K[x^p, \der^p]_{f^p} x^i\der^i.
\end{equation}
 In particular, $ Z(L)_{f^p}=Z(A_1)_{f^p}=K[x^p, \der^p]_{f^p}$, and  so we can write
\begin{eqnarray*}
\Spec (Z(L))\backslash V(f^p)&=&\Spec (Z(L)_{f^p})=\Spec (Z(A_1)_{f^p})=\Spec (K[x^p, \der^p]_{f^p})\\
&=&\Spec (Z(A_1))\backslash V(f^p)=\Spec (K[x^p, \der^p])\backslash V(f^p).
\end{eqnarray*} 

Theorem \ref{pKGLD} gives explicit descriptions of the sets of prime, completely prime, primitive and maximal  ideals of the algebra $\L $. Let $\Spec_c (\L , {\rm ht}=1)$ be the set  of  completely prime ideals  of height 1 of the algebra $\L (f)$.

 \begin{theorem}\label{pKGLD}
Let $K$ be a field of characteristic $p>0$, $\L =K[x][y; \d := f\frac{d}{dx} ]$ where $f\in K[x]\backslash K$. Let $f=p_1^{n_1}\cdots p_s^{n_s}$ be a unique (up to permutation) product of irreducible polynomials of $K[x]$. Then 
\begin{enumerate}
\item The elements $p_1, \ldots , p_s$ are regular normal elements of the algebra $\L$ (i.e. $p_i$ is a non-zero-divisor of $\L$ and $p_i\L = \L p_i$). 
\item $\min (f) = \{ (p_1), \ldots , (p_s)\}$. 
\item $\Spec_c (\L )=\{ 0 , \L p_i,  (p_i, q_i)\, | \, i=1, \ldots , s;\; q_i\in  {\rm Irr}_m(F_i[y])\}$  where $F_i:=K[x]/(p_i)$ is a field  and ${\rm Irr}_m(F_i[y])$ is the set of monic irreducible polynomials of  the polynomial algebra $F_i[y]$ over the field $F_i$ in the variable $y$. If, in addition, the field $K$ is an algebraically closed and $\l_1, \ldots , \l_s$ are the roots of the polynomial $f$ then $\Spec_c (\L )=\{ 0 , \L (x-\l_i),  (x-\l_i, y-\mu)\, | \, i=1, \ldots , s;\; \mu\in K\}$. 
\item $\Spec (\L)=\Spec_c(\L )\coprod \{ \L  \gp \, | \, \gp\in \Spec (Z (\L ))\backslash \{ (0), V(f^p) \}$.
\item For all $\gp \in \Spec (Z (\L ))\backslash V(f^p)$, $$k(\gp )\t_{Z(\L )}\L /(\gp )\simeq \bigoplus_{i,j=0}^{p-1}\k (\gp ) x^iy^i=\bigoplus_{i,j=0}^{p-1}\k (\gp ) x^i\der^i\simeq M_{p}(k(\gp )),$$ the algebra of $p\times p$ matrices over the field of  fractions $\k( \gp)$ of the domain $Z(\L )/\gp$. 
\item $\Max (\L )=\Prim (\L )=\{ (p_i, q_i), \L\gm \, | \, i=1, \ldots , s;\; q_i\in  {\rm Irr}_m(F_i[y]), \gm \in \Max (Z(\L ))\backslash V(f^p) \}$.
\item $\Spec_c (\L , {\rm ht}=1)=\{ (p_1), \ldots , (p_s)\}$. If, in addition $K=\bK$, then $\Spec_c (\L , {\rm ht}=1)=\{ (x-\l_1), \ldots , (x-\l_t)\}$ where $\{ \l_1, \ldots , \l_t\}$ is the set of roots of the polynomial $f$. 
\end{enumerate}
\end{theorem}

{\it Proof.} 1. Statements 1 is proven above.

2. Since 
\begin{equation}\label{LLpi}
\L/ \L p_i\simeq F_i[y]
\end{equation}
is a polynomial algebra with coefficients in the field $F_i$  (since $yx-xy=f\in \L p_i$), the ideal $\L p_i$ is a completely prime ideal of $\L$. Now, statement 2 follows from the equality of ideals $(f)=(p_1)^{n_1}\cdots (p_s)^{n_s}$.

5. Since $\gp \in \Spec (Z (\L ))\backslash V(f^p)$, the element $f^p$ is a unit in the field $k(\gp )$. Now, the first isomorphism and the equality in statement 5 follows from (\ref{Lffloc}). Then using the equalities, $$[y,x^i]=ix^{i-1}\;\; {\rm  and }\;\;  [y^i, x]=iy^{i-1},$$  and the fact that $k(\gp )$ is a field, we see that the algebra $\bigoplus_{i,j=0}^{p-1}\k (\gp ) x^i\der^i$ is a simple, central $k(\gp )$-algebra of dimension $p^2$ over the field $k(\gp )$. Therefore, it is isomorphic to the matrix algebra $M_p(k(\gp ))$.

3-4. The algebra $\L$ is a domain, hence $0\in \Spec_c (\L )$. We have seen in the proof of statement 2 that the ideals $(p_1), \ldots , (p_s)$ of the algebra $\L$  are completely prime ideals.  By (\ref{LLpi}),
$$V(f)=  \{  \L p_i,  (p_i, q_i)\, | \, i=1, \ldots , s;\; q_i\in  {\rm Irr}_m(F_i[y])\}\subseteq \Spec_c (\L).$$
Given a nonzero prime ideal $P$ of $\L$ such that $P\not \in V(f)=V(f^p)$. Then $P_f$ is a  nonzero prime ideal of the algebra $\L_{f^p}= A_{1,f^p}$. By Corollary  \ref{a16Mar20}, the intersection $\gp := P\cap Z(L)$ is a nonzero prime ideal of the centre $Z(\L )$  of  the algebra $\L$. By Theorem \ref{16Mar20}.(4), $\L = \bigoplus_{i,j=0}^{p-1}Z(\L )x^iy^i$. Now, by statement 5, $\L \gp\in \Spec (L)\backslash V(f)$ and the prime ideal $\L \gp$ is not completely prime. Now, statements 3 and 4 follows from statement 5. 

6. Statement 6 follows from statement 4.

7. Statement 7 follows from statement 3.  $\Box$\\

{\bf Classification of simple $\L (f)$-modules.}  
 For an $\L $-module $M$, we denote by $\ann_\L (M)$ the annihilator of the $\L$-module $M$. For an algebra $A$, we denote by $\widehat{A}$, the set of isomorphism classes of (left) simple $A$-modules. An isomorphism class of a simple $A$-modules $M$ is denoted by $[M]$.   Let elements $a_1, \ldots , a_n\in A$ be  generators for a left ideal $I$ of the algebra $A$. Then we write $I=A(a_1, \ldots , a_n)$.  Theorem \ref{17Mar20} is a classification of simple $\L (f)$-modules.

  \begin{theorem}\label{17Mar20}
Let $K$ be a field of characteristic $p>0$, $\L =K[x][y; \d := f\frac{d}{dx} ]$ where $f\in K[x]\backslash K$. Let $f=p_1^{n_1}\cdots p_s^{n_s}$ be a unique (up to permutation) product of irreducible polynomials of $K[x]$. Then 
\begin{enumerate}
\item the map
$$ \Max (L)\ra \widehat{\L }, \;\; \gM \mapsto L(\gM )$$
is a bijection with inverse $[M]\mapsto \ann_\L (M)$ where $L(\gM )$ is a unique (up to isomorphism) simple direct summand/submodule/factor module of the (simple) matrix algebra $\L / \gM $. In particular, for all $\gm \in \Max (\L)\backslash V(f^p)$, $\dim_K(L(\L \gm )) = p\cdot \dim_K(Z(\L ) / \gm )<\infty$. 

\item For each maximal ideal  $(p_i, q_i)$ of $\L (f)$, where  $i=1, \ldots , s$ and $q_i\in  {\rm Irr}_m(F_i[y])\}$, $$L(p_i, q_i)=\L / (p_i, q_i)\simeq K[y]/(q_i)$$ and $\dim_K(L(p_i, q_i))=\dim_K( K[y]/(q_i))=\deg_y(q_i)<\infty$. 

\item Suppose that the field $K= \bK$ is an algebraically closed.  For each maximal  ideal $ \L \gm$ of $\L (f)$, where $\gm \in \Max (Z(\L ))\backslash V(f^p)$, 
 \begin{eqnarray*}
 L(\L \gm )&\simeq& \L/\L (\gm , x-\sqrt[p]{\xi})= \bigoplus_{i,j=0}^{p-1}K y^i\bar{1} \\
&\simeq& A_{1,f^p}/A_{1,f^p}(\gm , x-\sqrt[p]{\xi})= \bigoplus_{i,j=0}^{p-1}K\der^i\hat{1}, 
\end{eqnarray*}
where $x^p-\xi \in\gm$ for a unique element $\xi \in K$, $\bar{1}=1+\L (\gm , x-\sqrt[p]{\xi})$ and $\hat{1}=1+A_{1,f^p}(\gm , x-\sqrt[p]{\xi})$, $\dim_K L(\L \gm )=p<\infty$. 
\end{enumerate}
\end{theorem}

{\it Proof.} 1. Statements 1 follows at once from Theorem \ref{pKGLD}.(5,6).

2. Statement 2 is obvious.

3.  Notice that $(x-\sqrt[p]{\xi}))^p=x^p-\xi \in \gm$. 

By Theorem \ref{16Mar20}.(4) and the choice of the ideal $\gm$,  we have that $$\L / \L \gm =\bigoplus_{i=0}^{p-1} Ky^i\t K[x]/(x^p-\xi) =\bigoplus_{i=0}^{p-1} Ky^i\t K[x]/((x-\sqrt[p]{\xi})^p),$$ direct sums of tensor products of vector spaces.  Hence, 
 $$ \L/\L (\gm , x-\sqrt[p]{\xi})\simeq \bigoplus_{i=0}^{p-1} Ky^i\t K[x]/(x-\sqrt[p]{\xi})
\simeq \bigoplus_{i,j=0}^{p-1}K y^i\bar{1} 
 $$
 is a $p$-dimensional $\L$-module that is annihilated by the maximal ideal $\gm$. By statement 1, it must be $L(\L \gm )$. Since $\gm \not\in V(f^p)$, the central element $f^p$ acts as a bijection on the module $L(\L \gm)$. Therefore, 
 $$L(\L \gm)=L(\L \gm)_{f^p}= \L_{f^p}/\L_{f^p} (\gm , x-\sqrt[p]{\xi})\simeq  A_{1,f^p}/A_{1,f^p}(\gm , x-\sqrt[p]{\xi})= \bigoplus_{i,j=0}^{p-1}K\der^i\hat{1}.\;\; \Box $$
 
 The action of the elements $x$ and $y$ on the $K$-basis $\{ y\bar{1}\, | \, i=0, \ldots , p-1\}$ of the $\L (f)$-module $L(\L \gm )$ of Theorem \ref{17Mar20}.(3) is given below:
 \begin{eqnarray*}
x \cdot \bar{1}& = & \xi^\frac{1}{p} \bar{1}, \\
x \cdot y^i\bar{1} &=& \xi^\frac{1}{p}y^i\bar{1}+\sum_{j=0}^{i-1} {i\choose j}\xi_{ij}y^j\bar{1}\;\; {\rm where}\;\; \xi_{ij} =(-1)^{i-j}\phi_{ij} (\xi^\frac{1}{p}), \;\; \phi_{ij}=\d^{i-j-1}(f)\in K[x], \\
 y \cdot y^i\bar{1}&=&  y^{i+1}\bar{1}\;\; {\rm where}\;\; 0\leq i\leq p-2, \\
y \cdot y^{p-1}\bar{1} &=& \rho \bar{1}\;\; {\rm where} \;\; y^p-\rho \in \gm \;\; {\rm for \; a\; unique \; element}\;\; \rho \in K.
\end{eqnarray*}

 {\bf The Krull and global dimensions of the algebra $\L (f)$.}
 
 \begin{theorem}\label{A16Mar20}
Let $K$ be a field of characteristic $p>0$, $\L =K[x][y; \d := f\frac{d}{dx} ]$ where $f\in K[x]\backslash K$. Let $f=p_1^{n_1}\cdots p_s^{n_s}$ be a unique (up to permutation) product of irreducible polynomials of $K[x]$. Then 
\begin{enumerate}
\item The Krull dimension of $\L$ is  $\Kdim (\L )=2$.
\item The global  dimension of $\L$ is  $\gldim(\L )=2$.
\end{enumerate}
\end{theorem}

{\it Proof.}  Let $\L = \L (f)$. 

1. By Theorem \ref{16Mar20}.(4), the algebra $\L $ is a finitely generated $Z(\L )$-module. Therefore, the Krull dimension of the algebra $\L $ is equal to the Krull dimension of the polynomial algebra $Z(\L )$ in two variables (Theorem \ref{16Mar20}.(3)), and statement 1 follows. \\

2.  By \cite[Theorem 7.5.3.(i)]{MR}, $\gldim (\L ) \leq \gldim (K[X])+1=1+1=2$. \\

 By (\ref{LLpi}), $\gldim (\L/ \L p_i)= \gldim ( F_i[Y])=1<\infty$. Now, by \cite[Theorem 7.3.5.(i)]{MR},
 $$\gldim (\L ) \geq \gldim (\L / \L p_i)+1\stackrel{(\ref{LLpi})}{=} \gldim ( F_i[Y])+1=1+1=2.$$
  Therefore, $\gldim  (\L )=2$. $\Box$


\section{Isomorphism problems and groups  of automorphisms for  Ore extensions $K[x][y; f\frac{d}{dx} ]$}\label{AUTR}

 In this section, a proof  Theorem \ref{p5Mar20} is given. It can be deduced from Theorem \ref{Alev-Dum-P3.6} but we give  a different proof. \\

 Let $K(x)$ be the field of rational functions in the variable $x$. Then the Ore extension  $B_1:= K(x)[\der ; \frac{d}{dx}]$ is the localization $B_1=(K[x]\backslash \{ 0\})^{-1}A_1$ of the Weyl algebra $A_1$ at the Ore set $K[x]\backslash \{ 0\}$. \\

Notice that 
\begin{equation}\label{SKGfK}
\mS (K)\rtimes G_f(K)=\{ \s_{\l , \mu , p}\, | \, \l \in K^\times, \mu \in K, p\in K[x]\}
\end{equation}
where 
$$\s_{\l , \mu , p} (x)=\l x+\mu \;\; {\rm and}\;\; \s_{\l , \mu , p} (y)=\l^{d-1}y+p $$
since $\s_{\l , \mu , p} =s_{\l^{-d+1}p}\s_{\l , \mu}$ where $d=\deg (f)$. \\

{\bf Proof of Theorem \ref{p5Mar20}.} Let $\s$ be an automorphism of the $K$-algebra $\L = \L (f)$. It can be uniquely extended to a $\bK$-automorphism, say $\s$,  of the algebra $\bK\t_K\L$. Let $\l_1, \ldots , \l_t$ be the roots of the polynomial $f$ in $\bK$. By Theorem \ref{pKGLD}.(7), the automorphism $\s$ permutes the set $$\Spec (\bK\t_K\L , {\rm ht}=1)=\{ (x-\l_1), \ldots , (x-\l_t)\}$$ of height 1 prime ideals of the algebra $\bK\t_K\L$ that are generated  
 by  regular normal elements $x-\l_1, \ldots , x-\l_t$ 
of the domain $\bK \t_K\L$ and the set $\bK^\times$ is the group of units of the algebra $\bK\t_K \L$.  So,  we must have that 
$$\s (x) = \l x+\mu $$ for some elements $\l\in \bK^\times$ and $\mu \in \bK$.   Since $K[x]=\L \cap \bK [x]$, we must have that $$\s (x)\in \s (\L ) \cap \s (\bK [x] )= \L \cap \bK [x] =K[x],$$ and so  $\l \in K^\times$ and $\mu \in K$.  So, the automorphism $\s$ respects the polynomial algebra $K[x]$. In particular, it respects the Ore set $S=K[x]\backslash \{ 0\}$ of the algebra $\L$. The automorphism $\s$ can be uniquely extended to  an automorphism of the algebra $B_1=S^{-1}\L$. Then $\s (\der ) = \l^{-1}\der +q$ for some element $q\in K(x)$. In particular, 
$$\s (y) = \s (f\der )= \s (f) (\l^{-1} \der +q)=\l^{-1} \frac{\s (f)}{f} y +p\;\; {\rm where}\;\; p:=\s (f) q\in K[x] $$
and $\s (f) = \g f$ for some element $\g \in K^\times$. Clearly, $\g = \l^d$ where $d=\deg (f)$ is the degree of the polynomial $f$ (since $\s (x) = \l x+\mu$). So,
$$\s (x)=\l x+\mu \;\; {\rm and}\;\; \s (y)=\l^{d-1}y+p, $$
i.e. $\s \in \mS (K)\rtimes G_f(K)$, as required. $\Box$. \\

By Theorem \ref{p5Mar20}, 
\begin{equation}\label{SKGfK1}
G_f(K)=\mS (K)\rtimes G_f(K)=\{ \s_{\l , \mu , p}\, | \, \l \in K^\times, \mu \in K, p\in K[x]\}
\end{equation}
where the multiplication and the inversion in the group $G_f(K)$ are given by the rule (where $d=\deg (f)$):
\begin{eqnarray*}
 \s_{\l_1 , \mu_1 , p_1} \s_{\l_2 , \mu_2 , p_2}&=& \s_{\l_1\l_2 , \l_2\mu_1+\mu_2 , \l_2^{d-1}p_1+p_2},\\
\s_{\l , \mu , p}^{-1}&=&\s_{\l^{-1} , -\l^{-1}\mu , -\l^{-d+1}p}.
\end{eqnarray*}

{\bf The algebra $B_1$ and its automorphism group.} The element $f$ is a regular normal element of $\L$  (i.e. $\L f=f\L )$ since 
$$fy=yf-f'f=(y-f')f \;\;{\rm where}\;\;f'=\frac{df}{dx}.$$
 It determines the $K$-automorphism $\o_f$ of the algebra $\L $:
$$fu=\o_f(u)f,\;\; u\in \L ,$$
$$\o_f:x\mapsto x,\;\;y\mapsto y-f'.$$
We denote by $\L_f$  and  $A_{1,f}$ the localizations of the algebras  $\L $ and  $A_1$ at the powers of the element $f$, i.e.
$$\L_f=S^{-1}_f\L \;\;{\rm and  }\;\;A_{1,f}=S^{-1}_fA_1\;\;{\rm where}\;\;
S_f=\{ f^i\}_{i\geq 0}.$$
By (\ref{(2.1)}), 
\begin{equation}\label{(2.3)}
\L \subset A_1\subset \L_f= A_{1, f}=K[x, f^{-1}][\der ;\frac{d}{dx}]\subset B_1.
\end{equation}
Recall that   $ \Aut_K(K[x])=\{ \s_{\l , \mu}\, | \, \l \in K^\times , \mu \in K\}$, $\s_{\l ,\mu}(x)=\l x+\mu$ and  

\begin{eqnarray*}
 \Aut_K(K(x))&=&\{ \s_M\, | \, M\in {\rm PGL}_2(K)\}\simeq {\rm PGL}_2(K), \; \s_M\ra M\; {\rm  where}\;  
\s_M(x)=\frac{ax+b}{cx+d}, \\
M&=&\begin{pmatrix}
a & b \\  c &d
 \end{pmatrix}\in  {\rm PGL}_2(K):={\rm GL}_2 (K)/K^\times E\; {\rm  and}\; E=\begin{pmatrix}
1 & 0 \\ 0  &1
 \end{pmatrix}.
\end{eqnarray*}
The maps 
\begin{eqnarray*}
\Aut_K(K[x]) &\ra & \Aut_K(A_1), \;\; \s_{\l , \mu}\mapsto \s_{\l , \mu} : x\mapsto \l x+\mu , \;\; \der\mapsto \l^{-1}\der, \\
\Aut_K(K(x)) &\ra & \Aut_K(B_1), \;\; \s_M\mapsto \s_M : x\mapsto \frac{ax+b}{cx+d}, \;\; \der \mapsto \frac{\s_M(f)}{f\s_M(x)'}\der,
\end{eqnarray*}
 are group monomorphisms where $g'=\frac{dg}{dx}$ for  $g\in K(x)$. We identify these groups with their images, i.e. 
 $$ \Aut_K(K[x])\subseteq \Aut_K(A_1)\;\; {\rm and}\;\; \Aut_K(K(x))\subseteq \Aut_K(B_1). $$ 
 The automorphism group  $\Aut_K(B_1)$ acts in the obvious way on the algebra $B_1$. Let $$\mS_1:= \St_{\Aut_K(B_1)}(x):=\{ \s\in \Aut_K(B_1)\, | \, \s (x)=x\},$$ the stabilizer of the element $x\in B_1$ in $\Aut_K(B_1)$. Clearly,   
$$\mS_1=\{ s_q :\, | \, q\in K(x)\}\simeq (K(x), +), \;\; s_q\mapsto q$$  where $s_q (x)=x$ and $s_q(\der )=\der +q$.

\begin{lemma}\label{a2Mar20}
\begin{enumerate}
\item $\Aut_K(B_1) = \mS_1\rtimes \Aut_K(K(x))=\{ \s_{M, q}:=s_q\s_M\, | \, M\in {\rm PGL}_2(K), q\in K(x)\}$.
\item  $\Aut_K(B_1, K[x]):=\{ \s \in \Aut_K(B_1)\, | \, \s (K[x])=K[x]\} = \mS_1\rtimes \Aut_K(K[x])=\{ \s_{l,\mu, q}\, | \, \l\in K^\times, \mu\in K, q\in K(x)\}$ where $ \s_{l,\mu, q}(x)=\l x+\mu$ and  $ \s_{\l,\mu, q}(\der )=\l^{-1}\der+q$.
\end{enumerate}
\end{lemma}

{\it Proof}. 1. Since  $\mS_1:= \St_{\Aut_K(B_1)}(x)$, we must have $\mS_1\cap \Aut_K(B_1)=\{ e\}$ and $\s \mS_1\s^{-1}\subseteq \mS_1$ for all automorphisms $\s\in \Aut_K(B_1)$. Hence,  $\Aut_K(B_1) \supseteq  \mS_1\rtimes \Aut_K(K(x))$.

To prove that the reverse inclusion holds we have to show that every  element $\s \in \Aut_K(K(x))$ belongs to the group $\mS_1\rtimes \Aut_K(K(x))$. The group of units  $K(x)^\times := K(x)\backslash \{ 0\}$ of the algebra $B_1$ is a $\s$-invariant set, i.e. $\s (K(x)^\times)=K(x)^\times$. Hence so is  the field $K(x)$. Let $\tau$ be the restriction of the automorphism $\s$ to the field $K(x)$. Then $\s_1:=\tau^{-1}\s\in \mS_1$, and so $\s = \tau \s_1\in \mS_1\rtimes \Aut_K(K(x))$, as required.

2. Statement 2 follows from statement 1. $\Box $\\

Below is a different proof of Theorem \ref{Alev-Dum-P3.6} is given.\\

{\bf Proof of Theorem \ref{Alev-Dum-P3.6}.} Let $\s : \L (f) \ra \L (g)$ be an isomorphism of the $K$-algebras. It can be uniquely extended to a $\bK$-isomorphism  $\s :\bK\t_K \L (f) \ra \bK\t_K\L (g)$. Let $\l_1, \ldots , \l_{s'}$ (resp., $\l_1', \ldots , \l_t'$)  be the roots of the polynomial $f$ (resp., $g$)  in $\bK$. By Theorem \ref{pKGLD}.(7), the automorphism $\s$ maps bijectively  the set $\{ (x-\l_1), \ldots , (x-\l_{s'})\}$ of height 1 completely prime ideals of the algebra $\bK\t_K \L (f)$ to the set $\{ (x-\l_1'), \ldots , (x-\l_t')\}$ of height 1 completely prime ideals of the algebra $\bK\t_K \L (g)$.  Therefore, $s'=t$.  Since the elements $x-\l_1', \ldots , x-\l_t'$ are regular normal elements of the domain $\bK \t_K\L (g)$ and the set $\bK^\times$ is the group of units of the algebra $\L (g) $, we must have that 
$$\s (x) = \l x+\mu $$ for some elements $\l\in \bK^\times$ and $\mu \in \bK$.   Since $K[x]=\L (g) \cap \bK [x]$, we must have that $\s (x)\in \s (\L (f)) \cap \s (\bK [x] )= \L (g) \cap \bK [x] =K[x]$, and so  $\l \in K^\times$ and $\mu \in K$.  So, the isomorphism $\s$ respects the polynomial algebra $K[x]$ of the algebras $\L (f)$ and $\L (g)$. In particular, it respects the Ore sets $S=K[x]\backslash \{ 0\}$ of the algebras $\L (f)$ and $\L (g)$, i.e. $\s (S)=S$. The isomorphism $\s$ can be uniquely extended to  an isomorphism $$\s : B_1=S^{-1}\L (f) \ra B_1=S^{-1}\L (g).$$ Then $\s (\der ) = \l^{-1}\der +q$ for some element $q\in K[x]$. In particular, 
$$\s (y) = \s (f\der )= \s (f) (\l^{-1} \der +q)=\l^{-1} \frac{\s (f)}{g} y +p\;\; {\rm where}\;\; p:=\s (f) q\in K[x] $$
and $\s (f) = \g g$ for some element $0\neq \g \in K[x]$. Applying the same argument for the isomorphism $\s^{-1}: \L (g) \ra \L (f)$, we have that $\s^{-1} (g) = \g_1 f$ for some element $0\neq \g_1 \in K[x]$.
Therefore, $$f=\s^{-1}\s (f)=\s^{-1}(\g g)=\s^{-1}(\g ) \g_1f,$$ and so $\g , \g_1\in K^\times$, and $\g_1=\g^{-1}$. 
Clearly,  $\g = \l^d$ where $d=\deg (f)$ is the degree of the polynomial $f$ (since $\s (x) = \l x+\mu$), and the theorem follows. Furthermore, 
$$\s (x)=\l x+\mu \;\; {\rm and}\;\; \s (y)=\l^{d-1}y+p. \;\; \Box $$


\section{The eigengroup $G_f$ of a polynomial $f$}\label{pAUTR}

For each nonscalar monic polynomial $f(x)$, Proposition \ref{A9Mar20}, 
Theorem \ref{A10Mar20}, Theorem \ref{A11Mar20}, Theorem \ref{B11Mar20} and Theorem \ref{C11Mar20} are explicit descriptions of the eigengroup $G_f(K)$ in the case when the field $K$ is algebraically closed.
The case of an arbitrary field is obtained from these results, Theorem \ref{A12Mar20}. The aim of this section is to prove these results.\\

{\bf The eigengroup $G_U(K)$.}\\

{\em Definition, \cite{Bav-AutOreExt}.}  Given a group,  a $G$-module $V$ over a field $K$ and a non-empty subset $U$ of $V$. The {\em eigengroup} of the set $U$ in $G$, denoted by $G_U(K)$, is the set of all elements of the group $G$ such that the elements of the set $U$ are eigenvectors of  them with eigenvalues in the  field $K$. Clearly, the eigengroup is a subgroup of $G$.\\

 Clearly, $$G_U=\bigcap_{u\in U} G_u$$ where $G_u:=G_{\{ u\} }=\{g\in G\, | \, gu=\l(g)u$ for some $\l (g)\in K\}$. If $K$ is a subfield of a field $K'$ then $G_U(K)\subseteq G_U(K')$ where $U$ is a subset of the $G$-module $K'\t_KV$ over the field $K'$. \\

{\bf Finite subgroups of $\Aut_K(K[x])$.} Let $K$ be a field of prime characteristic $p>0$, $\Fp =\Z /\Z p$ is the field that contains $p$ elements, for each $n\geq 1$, $ \F_{p^n}$ is the finite field that contains $p^n$ elements, $\FFp =\bigcup_{n\geq 1}\F_{p^n}$ is the algebra is closure of the field $\Fp$. Clearly, $\FFp\subseteq \bK$ and  group of roots of 1 in the field $\bK$ is $\FFp^\times :=\FFp\backslash \{ 0\}$. The group $\Aut_K(K[x])=\{ \s_{\l , \mu}\, | \, \l\in K^\times , \mu \in K\}$ where $\s_{\l , \mu } (x)= \l x+\mu$ and $$\Aut_K(K[x])\simeq \Sh(K)\rtimes \mT\simeq K\rtimes K^\times $$ where $\Sh (K) :=\{ \s_{1, \mu}\, | \, \mu \in K\}\simeq (K, +)$, $\s_{1, \mu } \mapsto \mu$ and $ \mT :=\{ \s_\l , 0\}\, | \, \l \in K^\times\} \simeq T^\times$, $\s_{\l , 0}\mapsto \l$ is the the algebraic 1-dimensional torus.

The set $ \mU =\mU (K)= K\cap \FFp^\times$ is the group of roots of 1 of the field $K$. The map $\mU\ra \mT$, $u\mapsto \s_{u, 0}$ is a group monomorphism and we identify the group $\mU$ with its image, i.e. $\mU = \{ \s_{u, 0}\, | \, u\in \mU\}$. Let ${\rm or} (g)$  be the order of an element of a group $G$.

\begin{lemma}\label{aA7Mar20}

The group $\Sh (K)\rtimes \mU (K)=\{ \s_{u, \mu}\, | \, u\in \mU (K),  \mu \in K\}$ is the set of all finite order automorphisms of the group  $\Aut_K(K[x])$. The order of the element  $\s_{u, \mu}=\s_{u, 0}\s_{1, \mu}$ is 
$$ {\rm or} (\s_{\l , \mu})=\begin{cases}
{\rm or} (\l )& \text{if } \l\neq 1,\\
p& \text{if } \l=1,\mu \neq 0,\\
1& \text{if }\l=1, \mu = 0.\\
\end{cases}
$$
\end{lemma}

{\it Proof}.  For all $\l \in K^\times\backslash \{ 1\}$ and $\mu \in K$, $\s_{1, \mu}^i=\s_{1, i\mu}$ and $\s_{\l , \mu}^i=\s_{\l^i, \frac{1-\l^i}{1-\l}\mu}$, and statement 1 follows.  $\Box $

 By Lemma \ref{aA7Mar20}, 
\begin{equation}\label{ShUK}
1\ra \Sh (K)\ra \Sh (K)\rtimes \mU (K) \stackrel{\v}{\ra}  \mU (K) \ra 1, \;\;{\rm where}\;\;  \v (\s_{\l , \mu})= \l ,
\end{equation}
 is a short exact sequence of group homomorphisms.

 \begin{lemma}\label{a8Mar20}
If an element  $\l\in \mU (\bK)$  is a primitive $n$'th root of unity then $\Fp (\l ) =\F_{p^m}$ where $m=\min \{ k\geq 1 \, | \; n|(p^{k}-1)\}=$ the degree of the minimal polynomial of $\l$ over the field $\Fp$. 
\end{lemma}

{\it Proof}. The element $\l$ is algebraic over the field $\Fp$. Let $\v$ be its minimal polynomial over $\Fp$. Then the field $\Fp (\l ) \simeq \Fp [x]/(\v )$ is a finite field, and so $\Fp (\l ) = \F_{p^m}$ for some $m\geq 1$. 
 Now, $$\deg (\v ) =[\Fp (\l ):\Fp]=[\F_{p^m} :\Fp]=m.$$ Notice that the order of the group $\langle \l \rangle$, which is $n$ (since $\l $ is a primitive $n$'th root of unity) divides the order of the group $\F_{p^m}^\times$, which is $p^m-1$. Clearly,  $m\geq m':=\min \{ k, | \, n|(p^{k}-1)\}$.  We claim that $m=m'$. Suppose that $m>m'$, wee seek a contradiction. Then $\l^{p^{m'}}=\l$, and so $\l \in \F_{p^{m'}}$, hence $F_{p^m} = \Fp (\l ) \subseteq \F_{p^{m'}}$. Therefore, $m| m'$, a contradiction.  $\Box $\\

The next theorem is a classification of all the finite subgroups of   the automorphism group $\Aut_K(K[x])$.

\begin{theorem}\label{C7Mar20}
Let $G$ be a finite subgroup of $\Aut_K(K[x])$. Then $G=\widetilde{G}\rtimes \oG$ where $\widetilde{G}=G\cap \Sh (K)=\{\s_{1, \mu}\, | \, \mu \in V\}$, $V\subseteq K$ is a finite dimensional $\Fp (\l_n)$-subspace of $K$ and $\oG=\langle \s_{\l_n, (1-\l_n)\nu}\rangle $ is a cyclic group of order $n$ where $\l_n$ is a  primitive $n$'th root of 1 and $\nu \in K$.  In particular, the order of the group $G$ is $np^l$ where $l=\dim_{\Fp}(V)$ such that $m|l$ where $\Fp (\l_n)=\F_{p^m}$ for some $m\geq 1$ (Lemma \ref{a8Mar20}).  Conversely, given a finite dimensional $\Fp (\l_n)$-subspace $V$ of $K$,  an automorphism $\s_{\l_n, (1-\l_n)\nu}$  where $\l_n\in K$ is a   primitive $n$'th root of unity and $\nu \in K$. Let $\widetilde {G}:= \{\s_{1, \mu}\, | \, \mu \in V\}$ and  $\oG:=\langle \s_{\l_n, (1-\l_n)\nu}\rangle $. Then the semidirect product $\widetilde{G}\rtimes \oG$ is a finite subgroup of  $\Aut_K(K[x])$ of order  $np^l$ where $l=\dim_{\Fp}(V)$ such that $m|l$. 
 The element $\nu \in K$ is unique up to adding an arbitrary element of $V$, i.e. $\widetilde{G}\rtimes \langle \s_{\l_n, (1-\l_n)\nu}\rangle  \simeq \widetilde{G}\rtimes \langle \s_{\l_n, (1-\l_n)\nu'}\rangle $ iff $\nu'-\nu \in V$. Furthermore, $G=\{ \s_{\l_n, (1-\l_n)\nu}^i\s_{1,v}\, | \, 0\leq i\leq n-1, v\in V\}$ and $\s_{\l_n, (1-\l_n)\nu}^i\s_{1,v}=\s_{\l_n^i, (1-\l_n^i)\nu}\s_{1,v}: x\mapsto \l_n^ix+(1-\l_n^i)\nu +v$.
\end{theorem}

{\it Proof}.  Let $G$ be a finite subgroup of $\Aut_K(K[x])$. Then, by (\ref{ShUK}), the group  $\v (G)$ is a finite subgroup of $\mU (K)$ of order $n$, hence $\v (G)=\langle \l_n\rangle$ where $\l_n$ is a  primitive $n$'th root of 1. Fix an element, say $\s_{\l_n, (1-\l_n )\nu}\in G$ where $\nu \in K$, such that $\v (\s_{\l_n, (1-\l_n )\nu})=\l_n$. Then $\oG=\langle \s_{\l_n, (1-\l_n)\nu}\rangle $ is a cyclic group of order $n=|\oG |$ since  
$\s_{\l_n, (1-\l_n)\nu}^i=\s_{\l_n^i, (1-\l_n^i)\nu}$ for all $i\geq 1$. Therefore, $$G=\widetilde{G}\rtimes \oG\;\; {\rm  where }\;\; \widetilde{G}:=G\cap \Sh (K)=\{\s_{1, \mu}\, | \, \mu \in V\},$$ $V\subseteq K$ is a finite dimensional $\Fp$-subspace of $K$ since $\s_{1, \mu}^i=\s_{1, i\mu}$ for all $i\geq 0$.  Furthermore, $\l_n V\subseteq V$, i.e. the $\Fp$-vector space $V$ is a $\Fp (\l)$-module  since 
$$\s_{\l_n, (1-\l_n)\nu}^{-1}\s_{1, \mu}\s_{\l_n, (1-\l_n)\nu}=\s_{1, \l\mu}.$$

 Clearly, $|G|= |\widetilde{G}| | \oG |=p^ln$ and $m|l$ since $V$ is a $\F_{p^m}$-module and $\dim_{\F_{p^m}}(V)=\frac{l}{m}$. 

The converse, is obvious. 

Clearly,  $\widetilde{G}\rtimes \langle \s_{\l_n, (1-\l_n)\nu}\rangle  \simeq \widetilde{G}\rtimes \langle \s_{\l_n, (1-\l_n)\nu'}\rangle $ iff there is natural number $i$ such that $1\leq i\leq n-1$ and a vector $v\in V$ such that 
$$\s_{\l_n, (1-\l_n)\nu'}=\s_{\l_n, (1-\l_n)\nu}^i\s_{1,v}=\s_{\l_n^i, (1-\l_n^i)(\nu +(1-\l_n^i)^{-1}v)}$$ 
iff $i=1$ and 
$\nu'=\nu +(1-\l_n^i)^{-1}v\in V$ iff $\nu' -\nu \in V$ since $(1-\l_n^i)^{-1}V=V$. $\Box $\\

The automorphism group $\Aut_K(K[x])$ acts on the set $\Max (K[x])$ of maximal ideals of $K[x]$ in the obvious way. If $K=\bK$ then $\Max (K[x])=\{ (x-\g )\, | \, \g \in K\}$ and the action takes the form: For all $\s \in \Aut_K(K[x])$ and $\g \in K$, 

\begin{equation}\label{sxgm}
\s ((x_\g ))=(x-\s^{-1} (\g ))\;\; {\rm where}\;\; \s^{-1} (\g ):=\s^{-1} (x)|_{x=\g }.
\end{equation}
Let us identify the set $\Max (K[x])$ with $K$ via $(x-\g )\mapsto \g$. Then the action of the group $\Aut_K(K[x])$ on $\Max (K[x])=K$ is given below:
\begin{equation}\label{sxgm1}
\Aut_K(K[x])\times K\ra K, \;\; (\s , \g)\mapsto \s*\g :=\s^{-1} (\g )=\s^{-1} (x)|_{x=\g }.
\end{equation}
If $\s = \s_{\l , \mu}$ then $\s_{\l , \mu}^{-1}=\s_{\l^{-1}, -\l^{-1}\mu}$ and $\s_{\l , \mu}*\g = \s_{\l , \mu }^{-1}(\g ) = \s_{\l^{-1}, -\l^{-1}\mu}(\g)=\l^{-1}\g -\l^{-1}\mu$. 

Every automorphism $\s_{\l \mu}\in \Aut_K(K[x])$ with $\l \neq 1$ can be uniquely written in the form $\s_{\l , (1-\l )\nu}$ where $\nu =(1-\l )^{-1}\mu$. Notice that 
$$ \s_{\l , (1-\l )\nu}*(\nu) = \nu .$$
Furthermore, the set $\{ \nu \}$ is the only 1-element orbit in $K$ of the cyclic group $\langle 
\s_{\l , (1-\l )\nu}\rangle$ generated by the automorphism $\s_{\l , (1-\l )\nu}$. The number of elements in any other orbit is equal to the order of the group  $\langle 
\s_{\l , (1-\l )\nu}\rangle$ which is the order of the element $\l$ in the group $(K^\times , \cdot )$. 

Suppose that $K=\bK$. Let $f\in K[x]$ be a nonscalar monic polynomial that has at least two distinct roots in $K=\bK$. Recall that $\CR (f)$ is the set of all roots of the polynomial $f$ counted with multiplicity and $\CR_d(f)$ be the set of all {\em distinct} roots of $f$ (i.e. each root has multiplicity 1). Example. For $f=(x-1)^2(x-2)^3$, $\CR (f)=\{ 1, 1, 2, 2, 2\}$ and $\CR (f)=\{ 1, 2\}$.

 The group $G_f$ permutes the roots in $\CR (f)$ and $\CR_d(f)$ via the action (\ref{sxgm}). Let us stress that the action of $G_f$ on $\CR (f)$ respects the multiplicity. If the group $G_f$ is finite then  $G_f=\widetilde{G}_f\rtimes \oG_f$ and $\widetilde{G}_f=\Sh_V$ is a normal subgroup of $G_f$. For a set $\CR =\CR (f), \CR_d(f)$ and a group $G=G_f, \widetilde{G }_f, \oG_f$, we denote by $\CR / G$ is the set of $G$-orbits in $\CR$. \\

{\bf Invariants and eigenalgebras of finite subgroups of $\Aut_K(K[x])$.}  Notice that 
\begin{equation}\label{fp=prod1}
x^p-x=\prod_{i\in \Fp}(x-i).
\end{equation}

For each element $\mu \in K^\times$, let 
\begin{equation}\label{fp=prod}
f_\mu (x):= \prod_{i=0}^{p-1}\s_{1, \mu }^i(x) = \prod_{i=0}^{p-1}(x-i\mu)=\prod_{i\in \Fp}(x-i\mu)=x^p-\mu^{p-1}x. 
\end{equation}
The equality above follows at once from (\ref{fp=prod1}): $f_\mu (x)=\prod_{i=0}^{p-1}(x-i\mu)=\mu^p\prod_{i=0}^{p-1}(\mu^{-1}x-i)=\mu^p((\mu^{-1}x)^p-\mu^{-1}x)=x^p-\mu^{p-1}x$. For all $\alpha , \beta \in \Fp$:
$$f_\mu (\alpha x+\beta x') = \alpha f_\mu(x) +\beta f_\mu (x')$$
since $\g^p=\g$ for all $\g \in \Fp$. In particular, the map $K\ra K$, $\l \mapsto f(\l)$ is a $\Fp$-linear map. Hence for all elements $\l\in K$, 
\begin{equation}\label{fp=prod2}
x^p-\mu^{p-1}x -(\l^p-\mu^{p-1}\l)=f_\mu (x)-f(\l) = f_\mu (x -\l ) =\prod_{i=0}^{p-1}(x-\l -i\mu)=\prod_{i\in \Fp}(x-\l -i\mu).
\end{equation}
Given a $K$-algebra $A$, and a subgroup $G$ of the automorphism group $\Aut_K(A)$. A group homomorphism $\chi : G\ra K^\times$ is called a {\em character} of the group $G$ in $K$. Let $\widehat{G} (K)$  be the (multiplicative) {\em group of characters} of the group $G$ in $K^\times$ . The multiplication in the group $\widehat{G}(K)$ is given by the rule:  For all  $\chi , \psi \in \widehat{G}(K)$, $(\chi \psi ) (g) = \chi (g) \psi (g)$ for all elements $g\in G$. \\

{\em Definition.} The direct sum of $G$-eigenspaces,  
$$\mE (A)=\mE (A, G):=\bigoplus_{\chi\in \widehat{G}(K)} A^\chi, \;\; {\rm where}\;\; A^\chi :=\{ a\in A\, | \, g(a) = \chi (g) a \;\; {\rm for\; all}\;\; g\in G\},$$
is called the $G$-{\em eigenalgebra} of $A$. The direct sum is a $\widehat{G}(K)$-graded algebra since $A^\chi A^\psi \subseteq A^{\chi \psi}$ for all $\chi , \psi \in \widehat{G}(K)$.\\

 If $e$ is the identity element of the character group $\widehat{G}(K)$ then $A^e=A^G$ is the {\em algebra of $G$-invariants}. In particular $A^G\subseteq \mE (A,G)$.  

The set $\Supp (A,G):=\{ \chi \in \widehat{G}(K)\, | \, A^\chi \neq 0\}$ is called the {\em support} of $G$ in $A$.  For each character $\chi\in \Supp (A, G)$, the vector space $A^\chi$ is called the $\chi$-{\em weight/eigenvalue subspace} 
of the algebra $A$. If the algebra $\mE (A, G)$ is a domain (eg, the algebra $A$ is a domain) then the support $\Supp (A, G)$ is a submonoid of  $\widehat{G}(K)$ and the algebra $\mE (A, G)$ is a $\Supp (A, G)$-graded algebra. 
If the algebra $A$ is a commutative algebra then the {\em Frobenius} endomorphism $\Fr : A\ra A$, $a\mapsto a^p$ is a $\Fp$-algebra endomorphism of $A$. It is a monomorphism if the algebra $A$ is a domain.  By the very definition,  the Frobenius endomorphism commute with all endomorphisms of the ring $A$.

\begin{lemma}\label{a9Mar20}
Let $A$ be a commutative $K$-algebra and $G$ be a subgroup of the automorphism group $\Aut_K(A)$. 
\begin{enumerate}
\item The algebras $\mE (A, G) $ and $A^G$ are $\Fr$-stable (that is $\Fr (\mE (A, G))\subseteq \mE (A, G) $ and $\Fr (A^G)\subseteq A^G$). 
\item Suppose that  the algebra $A$ is reduced and  $\Fr (K)=K$.  If $g(\Fr (a)) = \chi (g) \Fr (a)$ for all $g\in G$ then $g(a) = (\chi (g))^\frac{1}{p}a$. In particular, $\Fr \in \Aut_{\Fp}(\mE (A, G))$ and $\Fr \in \Aut_{\Fp}(A^G)$.
\end{enumerate}
\end{lemma}

{\it Proof}. 1. The Frobenius endomorphism commutes with all endomorphisms of the ring $A$, and statement 1 follows.

2. The equality $\Fr (K)=K$ implies that $\Fr\in \Aut_{\Fp}(K)$. Since the  algebra $A$ is reduced the Frobenius endomorphism $A$ is a monomorphism. If $g(\Fr (a)) = \chi (g) \Fr (a)$ for all $g\in G$ then  $(g(a) - (\chi (g))^\frac{1}{p}a)^p=0$, and so $g(a) = (\chi (g))^\frac{1}{p}a$ for all $g\in G$, and statement 2 follows.  $\Box $\\

Let $V\subseteq K$ be a  $\F_{p^m}$-subspace of $K$.  The subgroup  $\Sh_V:=\{ \s_{1, v}\, | \, v\in V \}$ of $\Aut_K(K[x])$ is called the {\em shift group} that is determined by the  $\F_{p^m}$-subspace $V$.  Proposition \ref{B7Mar20} describes the  algebra  of invariants and the eigenalgebra of  the shift group $\Sh_V$.

\begin{proposition}\label{B7Mar20}

Let $V\subseteq K$ be a nonzero $\F_{p^m}$-subspace of $K$.  
Then $\mE (K[x], \Sh_V)=K[x]^{\Sh_V}$.
\begin{enumerate}
\item If $\dim_{\F_{p^m}}(V)=\infty$ then $K[x]^{\Sh_V}=K$. 
\item If $l=\dim_{\F_{p^m}}(V)<\infty$ and $\{ \mu_1, \ldots , \mu_l\}$ is a basis  of the vector space $V$ over $\F_{p^m}$  then 
\begin{enumerate}
\item the fixed algebra $K[x]^{\Sh_V}=K[f_V]$ is a polynomial algebra in $f_V:=\prod_{v\in V} (x-v)$,  the polynomial $f_V$ is divisible by  the polynomial $\prod_{i=1}^l f_{\mu_i}$,  
\item for all elements $\alpha , \beta \in \F_{p^m}$ and $\l\in K$, 
$f_V(\alpha x+\beta \l ) = \alpha f_V(x)+\beta f_V(\l)$. In particular, the map $K\ra K$, $\l \mapsto f_V(\l)$ is a $\F_{p^m}$-linear map, 
\item If $V\subset V'$ are distinct  $\F_{p^m}$-subspaces       of $K$ then $f_V | f_{V'}$ and $f_V \neq  f_{V'}$, and 
\item $\frac{df_V}{dx}\neq 0$. 

\end{enumerate}
\item In particular, for all elements $\mu \in K^\times$, $K[x]^{\s_{1, \mu}}=K[f_\mu ]$. 
\end{enumerate}
\end{proposition}

{\it Proof}. For all elements $\mu \in K$, the map $$\s_{1, \mu}-1: K[x]\ra K[x], \;\; \psi (x) \mapsto \psi (x+\mu ) -\psi (x)$$ is locally nilpotent map. Therefore, the element 1 is the only eigenvalue for the map $\s_{1, \mu}$, and so 
$\mE (K[x], \Sh_V)=K[x]^{\Sh_V}$.

1. Statement 1 follows from statement 2. Let $\{ \mu_i\}_{i\in \N}$ be an $\F_{p^m}$-linearly independent elements of the vector space $V$ and $V_i=\bigoplus_{j=1}^i\F_{p^m} \mu_i$. Then $V_1\subset V_2\subset \cdots \subset 
V_\infty :=\bigoplus_{i\geq 1}\F_{p^m}\mu_i\subseteq V$. Hence, 
$$K[x]^{\Sh_{V_1}} \supseteq K[x]^{\Sh_{V_2}}\supseteq \cdots \supseteq 
K[x]^{\Sh_{V_\infty} }=\bigcap_{i\geq 1}K[x]^{\Sh_{V_i}} =\bigcap_{i\geq 1}K[f_{V_i}] = K \supseteq K[x]^{\Sh_{V} }\supseteq K, $$
and so $K[x]^{\Sh_V}=K$.

2.(a,b). For all elements  $v'\in V$, 
$$\s_{1, v'}(f_V)=\prod_{v\in V} (x-v'-v)=f_V.$$ Therefore, $K[x]^{\Sh_V}\supseteq K[f_V]$. By (\ref{fp=prod}), the polynomial $\prod_{i=1}^l \prod_{u\in \Fp} (x-u\mu_i) =\prod_{i=1}^l f_{\mu_i}$ is a divisor of the polynomial $f_V$. In particular, $f_V(0)=0$.  

First,  we prove that the statements (a) and (b) hold in the case when $K=\bK$ and then we deduce that the statements (a) and (b) hold for an arbitrary field $K$. 

So, suppose that  $K=\bK$.  Let $g(x)\in K[x]^{\Sh_V}$ be a noscalar monic polynomial. Let $\g $ be a root of $g(x)$. Then for all elements $v\in V$, the element $\g +v$ is also a root of the polynomial $g(x)$. So, the set of all roots of the polynomial $g(x) $ is a disjoin union of the sets $\coprod_{i=1}^s \{ \g_i+V\}$ for some roots $\g_i$ of $g(x)$. 

 Therefore, the polynomial 
$$g(x)=\prod_{i=1}^s \prod_{v\in V} (x-\g_i -v) =\prod_{i=1}^s f_V(x-\g_i )$$
is a product of $\Sh_V$-invariant polynomials $f_V(x-\g_i)$ of the same  degree $p^{lm}$. In particular, every nonscalar $\Sh_V$-invariant polynomials has degree at least $p^{lm}$. Therefore, for all elements $\l , \mu\in K$, 
the difference $$c_{\l , \mu}:=f_V(x+\l) -f_V(x+\mu)$$ of two monic $\Sh_V$-invariant polynomials of degree $p^{lm}$ must be a constant which is equal to 
$f_V(\l) -f_V(\mu)$.  Therefore, 
$$f_V(x+\l) -f_V(\l)=f_V(x+\mu) -f_V(\mu).$$
When $\mu =0$, we have that $$f_V(x+\l)=f_V(x) +f_V(\l) -f_V(0)=f_V(x) +f_V(\l)$$ since $f_V(0)=0$. Since for all elements $u\in \F_{p^m}^\times$, 
$$f_V(u x)=u^{p^{lm}}\prod_{v\in V}(x-u^{-1}v) = u\prod_{v\in V}(x-v) =uf_V(x),$$  and $f_V(0 x)=0=0f_V(x)$, we see that for all elements $\xi \in \F_{p^m}$, $f_V(\xi x) = \xi f_V(x)$. Now, the statement (b) follows. 

Now, the polynomial 
$$g(x)=\prod_{i=1}^s f_V(x-\g_i )=\prod_{i=1}^s (f_V(x)-f_V(\g_i ))\in K[f_V(x)]$$
and the statement (a) follows. 

Suppose that $K$ is not necessarily algebraically closed field and $g(x)\in K[x]^{\Sh_V}$ be a noscalar monic polynomial. Then $g(x) \in \bK [f_V(x)]$. 
Since the $\bK = K\oplus W$ for some $K$-subspace of $K$ and $f_V(x) \in K[x]$, we must have that $g(x) \in K[f_V(x)]$ since $$\bK [f_V(x)]=K [f_V(x)]\oplus \bigoplus_{i\geq 0}Wf_V(x)^i\subseteq K[x]\oplus\bigoplus_{i\geq 0}Wf_V(x)^i .$$
Now, the statements (a) and (b) hold for the field $K$.

(c) The statement (c) follows from the statement (a).

(d) WLOG we may assume that $K=\bK$. Suppose that $\frac{df_V}{dx}=0$. Then $f_V=g^p$ for some polynomial $g$. This is not possible as every root of $f$ has multiplicity 1. 

3. Statement 3 is a particular case of statement 2. $\Box$\\

Notice that for all natural numbers $m\geq 1$,  
\begin{equation}\label{mfp=prod1}
x^{p^m}-x=\prod_{i\in \F_{p^m}}(x-i).
\end{equation}

By (\ref{mfp=prod1}), for each element $\mu \in K^\times$, let 
\begin{equation}\label{mfp=prod}
f_{ p^m, \mu }(x):= f_{\F_{p^m}\mu}(x)= \prod_{i\in \F_{p^m}}(x-i\mu)=x^{p^m}-\mu^{p^m-1}x. 
\end{equation}
 For all $\alpha , \beta \in \F_{p^m}$:
 \begin{equation}\label{mfp=prod3}
f_{ p^m, \mu } (\alpha x+\beta x') = \alpha f_{ p^m, \mu }(x) +\beta f_{ p^m, \mu } (x')
\end{equation}
 since $\g^{p^m}=\g$ for all $\g \in \F_{p^m}$. In particular, the map $K\ra K$, $\l \mapsto f_{ p^m, \mu }(\l)$ is a $\F_{p^m}$-linear map. Hence, for all elements $\l\in K$, 

\begin{equation}\label{mfp=prod2}
x^{p^m}-\mu^{p^m-1}x -(\l^{p^m}-\mu^{p^m-1}\l)=f_{ p^m, \mu } (x)-f_{ p^m, \mu }(\l) = f_{ p^m, \mu } (x -\l ) =\prod_{i\in \F_{p^m}}(x-\l -i\mu).
\end{equation}

Theorem  \ref{8Mar20} describes the  algebra  of invariants and the eigenalgebra of  a `generic' finite subgroup of $\Aut_K(K[x])$.
 
\begin{theorem}\label{8Mar20}
Let $G=\widetilde{G}\rtimes \oG$ be a finite subgroup of $\Aut_K(K[x])$ (Theorem \ref{C7Mar20})   where $\oG:=\langle \s_{\l_n, (1-\l_n)\nu}\rangle $ and $\widetilde {G}:= \{\s_{1, \mu}\, | \, \mu \in V\}$,  $\l_n\in K$ is a primitive $n$'th   root of unity, $n\geq 2$ and $\nu \in K$,  $V$ is a nonzero  finite dimensional $\Fp (\l_n)$-subspace of $K$, and $\Fp (\l_n)=\F_{p^m}$ for some $m\geq 1$ (Lemma \ref{a8Mar20}). Then 
\begin{enumerate}
\item $\s_{\l_n, (1-\l_n)\nu} (f_V(x-\nu ))=\l_n f_V(x-\nu )$. 
\item $K[x]^G=K[f_V^n(x-\nu )]$ is a polynomial algebra in $f_V^n(x-\nu ):=\Big( f_V(x-\nu ) \Big)^n$. 
\item The $G$-eigenvalue subalgebra of $K[x]$ is $\mE (K[x], G)=\bigoplus_{i=0}^{n-1}
f_V^i(x-\nu )K[x]^G$, a direct sum of distinct  $G$-eigenspaces. 
\end{enumerate}
\end{theorem}

{\it Proof}. 1. Let $\s = \s_{\l_n, (1-\l_n)\nu}$. Suppose that $l=\dim_{\F_{p^m}}(V)$. Then $\deg (f_V(x))=p^{lm}$. 
Now,  by Proposition \ref{B7Mar20}.(2b), 
$$
\s (f_V(x-\nu )) = f_V(\s (x-\nu )) =f_V(\l_n(x-\nu ))=\l_n f_V(x-\nu )$$
since $\l_n\in K(\l_n ) = \F_{p^m}$.  

2 and 3. By Proposition \ref{B7Mar20}.(2b), the  $\widetilde{G}$-eigenvalue subalgebra of $K[x]$ is the fixed algebra $K[x]^{\widetilde{G}}=K[f_V(x)]=K[f_V(x-\nu)]$ since $$f_V(x)=f_V(x-\nu +\nu) = f_V(x-\nu) +f_V(\nu).$$ 
 By statement 1, $\s(f_V(x-\nu ))=\l_n f_V(x-\nu )$, hence   the  $\widetilde{G}$-eigenvalue subalgebra of $K[x]$, $K[f_V(x-\nu)]$, is $\s$-invariant, and statements 2 and 3  follow.   $\Box $\\

Proposition  \ref{A8Mar20} describes the  algebra  of invariants and the eigenalgebra of  the subgroup $G=\langle \s_{\l_n, (1-\l_n)\nu}\rangle$ of $\Aut_K(K[x])$ where  $\l_n\in K$ is a   primitive $n$'th  root of unity, $n\geq 2$ and $\nu \in K$. 

\begin{proposition}\label{A8Mar20}
Let $G=\langle \s_{\l_n, (1-\l_n)\nu}\rangle$ be a finite subgroup of $\Aut_K(K[x])$ where  $\l_n\in K$ is a   primitive $n$'th  root of unity, $n\geq 2$ and $\nu \in K$. Then 
\begin{enumerate}
\item $\s_{\l_n, (1-\l_n)\nu} (x-\nu )=\l_n (x-\nu )$. 
\item $K[x]^G=K[(x-\nu )^n]$ is a polynomial algebra in $(x-\nu )^n$. 
\item $\mE (K[x], G)=K[x]=\bigoplus_{i=0}^{n-1}
(x-\nu )^iK[x]^G$ is a direct sum of distinct  $G$-eigenspaces.

\end{enumerate}
\end{proposition}

{\it Proof}. 1. Statement 1 is obvious.

 2. Statement 2 follows from statement 1 and the fact that $\l_n$ is a primitive $n$'th root of unity. 
 
 3. Statement 3 follows from statement 3. $\Box $\\

{\bf The eigengroup $G_f(K)$ of a polynomial $f\in K[x]$ that has single root in $\bK$.} 
For an element $\nu\in K$, the subset
\begin{equation}\label{mTnuK}
\mT_\nu (K):=\{ \s_{\l , (1-\l )\nu}\, | \, \l \in K^\times\}
\end{equation}
of $\Aut_K(K[x])$ is a subgroup which is isomorphic to the algebraic torus $\mT =(K^\times , \cdot )$ via 
$$\mT_\nu (K)\ra \mT, \;\;\s_{\l , (1-\l )\nu}\mapsto \l$$  since $\s_{\l , (1-\l )\nu}\s_{\l' , (1-\l' )\nu}=\s_{\l\l' , (1-\l )\nu}$ for all $\l , \l'\in K^\times $. 

Suppose that a monic nonscalar polynomial $f\in \bK [x]$ of degree $d$ has single root, say $\nu \in \bK$. Then  $d=p^rd_1$ where  for unique natural numbers $r$ and $d_1$ such that $p\nmid  d_1$. Then 
\begin{equation}\label{fxdn}
f=(x-\nu )^d=(x-\nu )^{p^rd_1}=(x^{p^r}-\nu^{p^r})^{d_1} = x^{p^rd_1}-d_1\nu^{p^r} x^{p^r(d_1-1)}+\cdots .
\end{equation}
Therefore,  $f=(x-\nu )^d\in K[x]$ iff $\nu^{p^r}\in K$ (since $p\nmid d_1$). 

The next proposition  describes all the monic polynomial $f\in K[x]$ such that $G_f(K) = \mT_\nu (K)$ for some $\nu \in K$. Furthermore, it describes the group $G_f$ for all polynomial $f\in K[x]$ that has a single root in $\bK$. 
 
\begin{proposition}\label{A9Mar20}
\begin{enumerate}
\item Let $f(x)\in \bK [x]$ be a monic nonscalar polynomial of degree  $d$. Then $f(x)=(x-\nu )^d$ for some $\nu \in \bK$ iff $G_f(\bK ) = \mT_\nu (\bK)$. 
\item  Let $f(x)\in K [x]$ be a monic nonscalar polynomial of degree  $d$ that has a single root $\nu \in \bK$.  Then 
$G_f(K)=\begin{cases}
\mT_\nu (K)& \text{if }\nu \in K,\\
\{ e\} & \text{if }\nu \not\in K.\\
\end{cases}$
\item Let $f(x)\in K [x]$ be a monic nonscalar polynomial of degree  $d$ that has a single root $\nu \in \bK$.  Then  $f(x)=(x-\nu )^d$ for some $\nu \in K$ iff $G_f(K ) = \mT_\nu (K)$. 
\end{enumerate}
\end{proposition}

{\it Proof}. 1. $(\Rightarrow )$ Suppose that $f=(x-\nu )^d$. An automorphism $\s\in \Aut_K(K[x])$
belongs to the group $G_f(\bK )$ iff $\s (x-\nu) =\l (x-\nu )$ for some element $\l \in \bK^\times $ iff $\s=\s_{\l , (1-\l )\nu}$ iff $G_f=\mT_\nu (\bK )$. 

$(\Leftarrow )$ Suppose that $G_f=\mT_\nu (\bK )$. Then $K[x]=\bigoplus_{i\geq 0} K (x-\nu )^i$ is a direct sum of the eigenspaces of the group $\mT_\nu (\bK )$. Therefore, $f(x)=(x-\nu )^d$ for some $d\geq 1$. 

2. Clearly,  $G_f(K)=G_f(\bK )\cap \Aut_K(K[x])=\mT_\nu (\bK ) \cap \Aut_K(K[x])$. By statement 1,   $G_f(K)\neq \{ e\}$ iff $e\neq \s_{\l , (1-\l )\nu}\in \mT_\nu (\bK ) \cap \Aut_K(K[x])$ where $1\neq \l\in K^\times$ and (necessarily) $\nu \in K^\times$ iff  $G_f(K)= \mT_\nu (K )$.

3. Statement 3 follows from statement 2. $\Box $\\

{\bf The eigengroup $G_f(K)$ of a polynomial $f\in K[x]$ that has  at least two distinct roots   in $\bK$.}  For each nonscalar monic polynomial $f(x)$,  Theorem \ref{A10Mar20}, Theorem \ref{A11Mar20}, Theorem \ref{B11Mar20} and   Theorem \ref{C11Mar20}  (resp.,  Theorem \ref{A12Mar20}) are explicit descriptions of the eigengroup $G_f(K)$ in the case when the field $K$ is algebraically closed (resp., in general case). \\

Lemma \ref{a10Mar20} is an explicit description of the roots of the polynomials of the type $g(f^n_V(x-\nu ))$.

\begin{lemma}\label{a10Mar20}
Suppose that $\l_n\in K$ is a primitive $n$'th   root of unity, $n\geq 2$ and $\nu \in K$,  $V$ is a nonzero  finite dimensional $\Fp (\l_n)$-subspace of $K$, and $K(\l_n)=\F_{p^m}$ for some $m\geq 1$ (Lemma \ref{a8Mar20}). Then 
\begin{enumerate}
\item For all elements $\rho \in K$, 
$$ f^n_V(x-\nu )-f_V^n(\rho )=\prod_{i=0}^{n-1}\prod_{v\in V} (x-\nu - \l_n^i\rho -v).$$
\item Let $g(x)=\prod_{j=1}^k(x-\xi_j)\in \bK [x]$ where $\CR (g)=\{ \xi_1, \ldots , \xi_k\}$ is the set of roots of the polynomial $g(x)$ counted with multiplicity. Then $\xi_j=f_V^n(\rho_j)$ for some element $\rho_j\in \bK$ and 
$$ g(f^n_V(x-\nu ))=\prod_{j=1}^k\prod_{i=0}^{n-1}\prod_{v\in V} (x-\nu - \l_n^i\rho_j -v).$$

\end{enumerate}
\end{lemma}

{\it Proof}. 1. By Proposition \ref{B7Mar20}.(2b),  the map $f_V$ is a $K(\l_n)$-linear map (since $ K(\l_n)=\F_{p^m}$), and the result follows:
\begin{eqnarray*}
 f^n_V(x-\nu )-f_V^n(\rho ) &=&  \prod_{i=0}^{n-1}(f_V(x-\nu )-\l_n^if_V(\rho ))= \prod_{i=0}^{n-1}f_V(x-\nu -\l_n^i\rho )\\
  &=& \prod_{i=0}^{n-1}\prod_{v\in V} (x-\nu - \l_n^i\rho -v).
\end{eqnarray*}
2. Notice that  $g(x)=\prod_{j=1}^k (f^n_V(x-\nu )-f_V^n(\rho_j ))$, and statement 2 follows from statement 1. $\Box$

\begin{theorem}\label{10Mar20}
Suppose that  a monic  polynomial   $f(x)\in K[x]$ has at least two distinct roots in $\bK$. Then the group $G_f(K)$ is a finite group, $G_f(K)=\widetilde{G}_f(K)\rtimes \oG_f(K)$ where $\oG_f(K)=\langle \s_{\l_n, (1-\l_n)\nu}\rangle$ and $\widetilde{G}_f(K)=\Sh_V(K)$, $\l_n\in K$ is a primitive $n$'th root of unity, $\nu \in K$, $V$ is a finite dimensional  $\mF_p(\l_n)$-subspace of $K$, and $\mF_p(\l_n)=\mF_{p^m}$ for some $m\geq 1$ (Lemma \ref{a8Mar20}). 
\end{theorem}

{\it Proof}. Since $G_f (K) = G_f(\bK ) \cap \Aut_K(K[x])$, it suffices to prove the theorem in the case when the field  $K$ is an algebraically closed field. So, we assume that $K=\bK$.  Let $\CR'$ be the set of  {\em distinct} roots of the polynomial $f$. The subgroup  $\widetilde{G_f}=G_f\cap \Sh (K)$ is equal to a group $ \Sh_V$ where $V$ is a finite dimensional vector space over the field $\Fp$ since $|V|\leq |\CR' |$ (as $\CR' +V\subseteq \CR' $). 

If $G$ is a finite subgroup of $G_f$ that contains the group $\Sh_V$ then $G=\Sh_V\rtimes \langle \s_{\l_n, (1-\l_n)\nu} \rangle$ where $\l_n$ is a primitive $n$'th root of unity and $\nu \in K$. The cyclic group $\langle \s_{\l_n, (1-\l_n)\nu} \rangle$ of order $n$ acts on the field $K$ and on the set $\CR'$, see (\ref{sxgm}). The point $\nu$ is the only  fixed point of the action and the orbit of every element $\l\neq \nu$ contains precisely $n$ elements. The polynomial $f$ contains at least two distinct roots. 
 Therefore $n\leq |\CR' |$. Then, 
by (\ref{ShUK}), the group $\v (G_f)$ is equal to  $\langle \s_{\l_{n'}, (1-\l_{n'})\nu'} \rangle$ where $\l_{n'}$ is a primitive $n'$'th root of unity and $\nu' \in K$. By Theorem \ref{C7Mar20}, $G_f=\Sh_V\rtimes \langle \s_{\l_{n'}, (1-\l_{n'})\nu'} \rangle$  is a finite group. $\Box $\\

Corollary \ref{b10Mar20} is a criterion for the group $G_f$ to be an infinite group. 

\begin{corollary}\label{b10Mar20}
Let $f(x)\in K[x]$ be  a nonscalar monic  polynomial. Then the group $G_f$ is an infinite group iff   $f(x)=(x-\nu )^d$ for some $\nu \in K$ iff $G_f(K ) = \mT_\nu (K)$ (see (\ref{mTnuK}) for the definition of the group $\mT_\nu (K)$).
\end{corollary}

{\it Proof}. The corollary follows from  Proposition \ref{A9Mar20}.(3)
 and Theorem \ref{10Mar20}. $\Box $\\
 
Recall that  the field $K$ is an algebraically closed field,   $f(x)\in K[x]$ is a monic  polynomial  that  has at least two distinct roots, and $\CR(f)$ is the set of all roots of $f$ counted with multiplicity.  Recall that (Theorem \ref{10Mar20} and Theorem \ref{C7Mar20}) the group $G_f=\widetilde{G_f}\rtimes \oG_f$ is a finite subgroup of $\Aut_K(K[x])$    where $\oG_f:=\langle \s_{\l_n, (1-\l_n)\nu}\rangle $ provided $\oG_f\neq \{ e\}$ and $\widetilde {G_f}:= \{\s_{1, \mu}\, | \, \mu \in V\}$,  $\l_n\in K$ is a primitive $n$'th   root of unity,  $\nu \in K$,  $V$ is a  finite dimensional $\Fp (\l_n)$-subspace of $K$, and $\Fp (\l_n)=\F_{p^m}$ for some $m\geq 1$ (Lemma \ref{a8Mar20}). \\

There are four distinguish cases:
\begin{enumerate}
\item $\widetilde{G_f}\neq \{ e\}$, $\oG_f\neq \{ e\}$,
\item $\widetilde{G_f}\neq \{ e\}$, $\oG_f=\{ e\}$,
\item $\widetilde{G_f}=\{ e\}$, $\oG_f\neq \{ e\}$,
\item $\widetilde{G_f}= \{ e\}$, $\oG_f=\{ e\}$.
\end{enumerate}
Below, in the case of $K=\bK$, for the polynomial $f$ criteria are given in terms of its roots for each case to hold. \\

{\bf A description of the group $\widetilde{G_f}(K)$ and a criterion for $\widetilde{G_f}(K)\neq \{ e\}$.} 

{\it Definition.} Let $f(x)\in K[x]$ be a nonscalar polynomial. Two distinct roots $\l , \l'\in \bK$ of the polynomial $f(x)$ are called a $K$-{\em shift pair} of $f(x)$ if 
\begin{equation}\label{FpllR}
\l-\l'\in K \;\; {\rm and}\;\; \Fp (\l -\l') +\CR (f) \subseteq \CR (f)
\end{equation}
where $\CR (f)$ is the set of all roots in $\bK$ of the polynomial $F(x)$ counted with multiplicity. \\

The set of all $K$-shift pairs of the polynomial $f(x)$ is denoted by ${\rm SP} (f, K)$. The vector space over $K$, 
\begin{equation}\label{FpllR1}
V(f, K)=\begin{cases}
\sum_{\{ \l , \l'\} \in {\rm SP} (f, K)}\Fp (\l -\l') & \text{if }{\rm SP} (f, K)\neq \emptyset,\\
0& \text{if }{\rm SP} (f, K)=\emptyset.\\
\end{cases}
\end{equation}
is called the $K$-{\em  shift vector space} of the polynomial $f(x)$.

For fields $K\subseteq L\subseteq \bK$, we have  
 ${\rm SP} (f, K)\subseteq {\rm SP} (f, L)\subseteq {\rm SP} (f, \bK)$ and $V(f, K)\subseteq V(f, L ) \subseteq V(f, \bK )$. 
 
 Proposition \ref{A25Mar20} gives an explicit description of the group $\widetilde{G}_f(K)$ and 
 a criterion for $\widetilde{G}_f(K)\neq \{ e\}$.
 
 \begin{proposition}\label{A25Mar20}
Let $f(x)\in K[x]$ be a nonscalar monic polynomial. Then 
\begin{enumerate}
\item  $\widetilde{G}_f(K)=\Sh_V$ where $V=V(f, K)$ is the $K$-shift vector space of $f$. 
\item $\widetilde{G}_f(K)=\{ e\}$ iff ${\rm SP}(f,K)=\emptyset$ iff for all distinct roots $\l, \l'\in \bK$ of the polynomial $f$ such that $\l - \l '\in K$ (if they exist), $\Fp (\l -\l') +\CR (f) \not\subseteq \CR (f)$. 
\item $\widetilde{G}_f(\bK )=\{ e\}$ iff ${\rm SP}(f,\bK )=\emptyset$ iff for all distinct roots $\l, \l'\in \bK$ of the polynomial $f$, $\Fp (\l -\l') +\CR (f) \not\subseteq \CR (f)$. 
\end{enumerate}
\end{proposition}

{\it Proof}. 1. It follows from the description of the group $\widetilde{G}_f(K)$ as a shift group,  $\widetilde{G}_f(K)=\Sh_V$,  that $V=V(f, K)$. 

2 and 3. Statement 2 follows from statement 1 and statement 3 is a particular case of statement 2.  $\Box $\\

{\it Definition.} Let $V$ be a nonzero finite dimensional $\Fp$-subspace of $K$. The largest finite field, denoted $\mF_{p^e}$, where $e=e(V)\geq 1$, such that $\mF_{p^e}V\subseteq V$ is called the {\em multiplier field} of $V$. The natural number $e=e(V)$ is called the $p$-{\em exponent} of $V$. \\

The multiplier  field $\mF_{p^e}$ is the composite of all finite fields $\mF_{p^m}$ such  that $\mF_{p^m}V\subseteq V$. If $\dim_{\Fp} (V)=n$
 then $p^m\leq |V|=p^n$, and so $m\leq n$.
  
 Let  $\l_{p^e-1}$ be  a primitive $p^e-1$'st root of unity (a generator of the cyclic group $\mF_{p^e}^\times$). Since $\l_{p^e-1}\in \mF_{p^e}$ and $\mF_{p^e}V\subseteq V$, we have that $\l_{p^e-1}V\subseteq V$.

 For each finite dimensional $\Fp$-subspace $V$ of the field $K$, Lemma \ref{b26Mar20} describes all the roots of unity $\l_n$ such that $\l_nV\subseteq V$.
 
 \begin{lemma}\label{b26Mar20}

Let $V$ be a nonzero finite dimensional $\Fp$-subspace of the field $K$, $\mF_{p^e}$ be its multiplier field. 
\begin{enumerate}
\item Suppose that  $\l_n$ is a primitive $n$'th root of unity. Then  $\l_nV\subseteq V$  iff  $n|p^e-1$.
\item $|\mF_{p^e}^\times |=1$ iff $(p,e)=(2,1)$ iff $\l_nV\subseteq V$ (where $\l_n$ is a primitive $n$'th root of unity) implies $\l_n=1$.
\end{enumerate}
\end{lemma}

{\it Proof}. 1. $\l_nV\subseteq V$ iff $\l_n\in \mF_{p^e}^\times$ iff  $n| |\mF_{p^e}^\times |$ iff  $n|p^e-1$.

2. Statement 2 follows from statement 1. $\Box $\\

{\bf Classification of  subgroups $G$ of 
 $\Aut_K(K[x])$ maximal  satisfying the property $G\cap \Sh (K)=\Sh_V$.}

 \begin{corollary}\label{a27Mar20}

Let $V$ be a nonzero finite dimensional $\Fp$-subspace of the field $K$, $\mF_{p^e}$ be its multiplier field and $\l_{p^e-1}$ be a primitive $p^e-1$'st root of unity. Then the finite groups $G_{V, \nu}:=\Sh_V\rtimes \langle \s_{\l_{p^e-1}, (1-\l_{p^e-1})\nu}\rangle$, where $\nu \in K/V$, are the maximal subgroups $G$ of the group $\Aut_K(K[x])$ that satisfy the property that $G\cap \Sh (K)=\Sh_V$.
\end{corollary}

{\it Proof}. Recall that for all elements $\l\in K^\times$ and  $\mu , v\in K$, $\s_{\l, \mu }\s _{1, v}\s^{-1}_{\l , \mu }= \s_{1, \l^{-1}v}$. Hence the groups $G_{V, \nu}$ are  well defined and every subgroup $H$ of $\Aut_K(K[x])$ such  that $H\cap \Sh (K)=\Sh_V$ is a finite group.

Given a finite subgroup $G'$ of $\Aut_K(K[x])$ such that $G'\cap \Sh (K)=\Sh_V$. By Theorem \ref{C7Mar20}, $G'=\Sh_V\rtimes \langle \s_{\l_n, (1-\l_n)\nu}\rangle$ for some $\nu \in K$ and a primitive $n$'th root of unity $\l_n$ such that $n| p^e-1$, by Lemma \ref{b26Mar20}.(1), and so $G'\subseteq G_{V, \nu}$. Since the groups $\{ G_{V, \nu }\}_{\nu \in K}$ are distinct, the corollary follows ($G_{V, \nu}=G_{V, \nu'}$ iff $\s_{\l , (1-\l )\nu'}=\s^i_{\l , (1-\l )\nu}\s_{1,v}$ for some natural number $i$ such that $1\leq i<p$ and $\gcd (i,p)=1$ and an element $v\in V$ where $\l =\l_{p^e-1}$ iff $\nu'=\nu+(1-\l^i)^{-1}v$ since $\s^i_{\l , (1-\l )\nu}\s_{1,v}=\s_{\l^i, (1-\l^i)(\nu+(1-\l^i)^{-1}v)}$  iff $\nu'\equiv \nu\mod V$ since $\Fp (\l) V= \mF_{p^e}V=V$). 
 $\Box $\\
 
 {\bf Criterion for $\widetilde{G}_f\neq \{e\}$ and $\oG_f\neq \{ e\}$.}
 
 \begin{lemma}\label{a28Mar20}

Suppose that $K$ is an algebraically closed field, $f(x)\in K[x]$ is monic nonscalar polynomial that has at least two distinct roots, $\widetilde{G}_f=\Sh_V\neq \{ e\}$ and  $\oG_f=\langle \s_{\l_n, (1-\l_n)\nu}\rangle \neq \{ e\}$ where $V$ is a nonzero $\Fp (\l_n)$-subspace of the field $K$ and $\l_n$ is primitive  $n$'th root of unity. Then 
\begin{enumerate}
\item $\l_n\in \mF_{p^e}$ where $\mF_{p^e}$ is the multiplier field of the $\Fp$-subspace $V$ of $K$, or, equivalently, $n|p^e-1$.
\item The group $G_f=\Sh_V\rtimes \oG_f$ is a subgroup of $\Sh_V\rtimes \langle \s_{\l, (1-\l)\nu}\rangle$ where $\l$ is a cyclic generator of the group $\mF_{p^e}^\times $, i.e. $\l = \l_{p^e-1}$ is primitive $p^e-1$'st root of unity. 
\end{enumerate}
\end{lemma}

{\it Proof}. 1. Since $\l_nV\subseteq V$, $\l_n\in \mF_{p^e}$ and the lemma follows from Corollary  \ref{a27Mar20}. $\Box $\\ 
 
{\em Definition.} Let $f(x) = x^d+a_{d-1}x^{d-1}+\cdots+a_1x+a_0\in \bK [x]$ be a monic polynomial of degree $d\geq 1$ where $a_i\in \bK$ are the coefficients of the polynomial $f(x)$. Then the natural number
$$ \gcd (f(x)):=\gcd \{ i\geq 1\, | \, a_i\neq 0\}$$
is called the {\em exponent} of $f(x)$.

Clearly, the exponent of $f(x)$ is the largest natural number $m\geq 0$ such that $f(x)=g(x^m)$ for some polynomial $g(x)\in K[x]$. \\

{\it Definition.} For a nonscalar polynomial $f\in K[x]$, we have the  unique product
\begin{equation}\label{gcdpf}
 \gcd (f) = p^s  \gcd_p (f)\;\; {\rm  where }\;\; s\geq 0, \;\; \gcd_p (f)\in \N\;\; {\rm  and}\;\;  p\nmid 
\gcd_p (f).
\end{equation}







\begin{proposition}\label{A26Mar20}

Suppose that $f(x) = f_V^i(x-\nu)$ for some nonzero finite dimensional $\Fp$-subspace $V$ of $K$, $\nu \in K$ and a natural number $i\geq 1$ (i.e. $\CR_d(f)=\nu +V$ and each root of $f(x)$ has multiplicity $i$). Let $\mF_{p^e}$ be the multiplier field of $V$, $\mF_{p^e}^\times =\langle \l_n\rangle$ where $n=p^e-1$. Then 
$$G_f=\begin{cases}
\Sh_V\rtimes \langle \s _{\l_n, (1-\l_n)\nu}\rangle\neq \Sh_V& \text{if }(p,e)\neq (2,1),\\
\Sh_V& \text{if }(p,e)= (2,1).\\
\end{cases}
$$
\end{proposition}

{\it Proof}. Since $|\CR_d(f)|=|\nu +V|=|V|\geq 2$,   the group $G_f=\widetilde{G}_f\rtimes \oG_f$ is a finite group. Clearly, $\Sh_V\subseteq \widetilde{G}_f$. In fact, $\Sh_V= \widetilde{G}_f$ since $\CR_d(f)=\nu +V$. 

 Suppose that $(p,e)\neq  (2,1)$. Recall that $n=p^e-1>1$ and $\mF_{p^e}^\times =\langle \l_n\rangle $, the multiplier field of $V$. In particular, $\l_nV\subseteq V$. Therefore, $\oG_f =\langle \s _{\l_n, (1-\l_n)\nu}\rangle$, by Lemma \ref{b26Mar20}.(1). 
 
 The case $(p,e)= (2,1)$ is obvious, see Lemma \ref{b26Mar20}.(2). $\Box $

 \begin{lemma}\label{a29Mar20}
Suppose that $\Fr (K)=K$ (where $\Fr (a)=a^p$, the Frobenius endomorphism). Then $G_{f^{p^n}}(K)=G_f(K)$ for all polynomials $f\in K[x]$ and all natural numbers $n\geq 1$. 
\end{lemma}

{\it Proof}. (i) $G_{f^{p^n}}(K)\supseteq G_f(K)$: If $\s \in G_f(K)$ then $\s (f) = \l f$ for some $\l \in K$, and so $\s (f^{p^n})=\l^{p^n}f^{p^n}$. This means that $\s \in G_{f^{p^n}}(K)$.

(ii) $G_{f^{p^n}}(K)\subseteq G_f(K)$:  If $\tau  \in G_{f^{p^n}}(K)$ then $\tau (f^{p^n}) = \mu f^{p^n}$ for some $\mu \in K$, and so $(\tau (f) -\mu^\frac{1}{p^n} f)^{p^n}=0$, i.e. $\tau (f) =\mu^\frac{1}{p^n} f$. This means that $\tau \in G_f(K)$. $\Box $\\

Let $f(x)=\sum a_ix^i\in K[x]$ be a monic nonscalar polynomial and $\gcd (f) =p^s\gcd_p(f)$. Suppose that $K=\bK$. Then there is a unique monic nonscalar polynomial $f_1(x)\in K[x]$ such that 
\begin{equation}\label{fxf1ps}
f(x)=f_1^{p^s}(x).
\end{equation}
Clearly, $\gcd (f_1) =\gcd_p(f)$, $\deg (f) = p^s\deg (f_1)$ and $f_1'\neq 0$ (the derivative of $f_1$).  \\

Theorem \ref{A10Mar20} is a criterion for the group $G_f=\widetilde{G_f}\rtimes \oG_f$ to have nontrivial subgroups $\widetilde{G_f}$ and $\oG_f$, it also gives an explicit description of the group $G_f$. 

\begin{theorem}\label{A10Mar20}
Suppose the field $K$ is an algebraically closed field and    $f(x)\in K[x]$ is a monic  polynomial  that  has at least two distinct roots,  $gcd (f) = p^s\gcd_p(f)$ and $f(x) =f_1^{p^s}(x)$ for a unique monic nonscalar polynomial $f_1(x)\in K[x]$, see (\ref{fxf1ps}). Suppose that  $V\neq 0$ is a  finite dimensional $\Fp $-subspace of $K$ and $\mF_{p^e}$ is the multiplier field of $V$. 
Then the following statements are equivalent:  
\begin{enumerate}
\item $\widetilde {G_f}=\Sh_V\neq \{ e\}$ and $\oG_f=\langle \s_{\l_n, (1-\l_n)\nu}\rangle \neq \{ e\}$.
\item  There is a primitive $n$'th root of unity $\l_n\neq 1$ 
 (in particular, $n\geq 2$) such that $\l_nV\subseteq V$, and an element $\nu \in K$ such that either $f(x)=f_V^i(x-\nu )$ for some natural number $i\geq 1$ and $n=p^e-1$ (in this case, $(p,e)\neq (2,1)$, $\widetilde{G}_f=\Sh_V$ and $\oG_f=\langle \s_{\l_{p^e-1}, (1-\l_{p^e-1})\nu}\rangle$ where $\mF_{p^e}^\times=\langle \l_{p^e-1}\rangle$) or otherwise $f(x)=f_V^i(x-\nu ) g(f_V^n(x-\nu )) $ for some natural number $i\geq 0$  and a  monic nonscalar polynomial $g(x)\in K[x]$ such that $g(0)\neq 0$, and the following two conditions hold:
 
\begin{enumerate}
\item $n\geq 2$ and $\gcd( \frac{p^e-1}{n},\gcd_p (g))=1$, and 
\item  $\CR (f) +\Fp (\l_n)(\l -\l')\not\subseteq \CR (f)$ for all distinct roots $\l$ and $\l'$ of the polynomial $f$  such that $\l -\l'\not\in V$.
\end{enumerate}
\end{enumerate}
Suppose that statement 1 holds. Then the natural number $i$ and the polynomial $g(x)$ in statement 2 are unique, the element $\nu$ is unique up to adding an arbitrary element of $V$ (i.e.  $\nu$ can be replaced by $\nu+v$ for any element $v\in V$), and the equality $f(x)=f_V^i(x-\nu ) g(f_V^n(x-\nu )) $ is unique (since $f_V(x-\nu -v)=f_V(x-\nu)$ for all $v\in V$). Furthermore,   $\s_{\l_n, (1-\l_n)\nu }(f)=\l_n^i f$, and  $f\in K[x]^{G_f}$ iff $n|i$.

In the second case, i.e.  $f(x)=f_V^i(x-\nu ) g(f_V^n(x-\nu )) $, $$n=\gcd (p^e-1, \gcd_p(h)$$ where $h(x)\in K[x]$ is a unique polynomial such that $f(x)=f_V^i(x-\nu ) h(f_V(x-\nu )) $ (i.e. $h(x) = g(x^n)$), and  either $\nu$ is a root of $f(x)$ (i.e. $i\geq 1$) or otherwise (i.e. $i=0$) $\nu$ is a root of $f_1'(x)$ (the derivative of $f_1(x)$).  
\end{theorem}

{\it Proof}. $(1\Rightarrow 2)$ Suppose that statement 1 holds. By Theorem \ref{10Mar20}, $\widetilde{G}_f=\Sh_V$ and $\oG_f=\langle \s_{\l_n, (1-\l_n)\nu}\rangle$ for a nonzero $\Fp$-subspace $V$ of $K$, a primitive $n$'th root of unity $\l_n\neq 1$ (in particular, $n\geq 2$) such that $\l_nV\subseteq V$ and an element $\nu \in K$. Then, by Theorem \ref{8Mar20}.(3), either  $f(x) = f^i_V(x-\nu )$ for some natural number $i\geq 1$ or otherwise 
$$f(x)=f_V^i(x-\nu ) g(f_V^n(x-\nu )) $$ for some natural number $i\geq 0$, $n\geq 2$,
 and a  monic nonscalar polynomial $g(x)\in K[x]$ such that $g(0)\neq 0$.

  In the first case, by Proposition \ref{A26Mar20}, $\widetilde{G}_f=\Sh_V\neq \{ e\}$ and $\oG_f=\langle \s_{\l_{p^e-1}, (1-\l_{p^e-1})\nu}\rangle\neq \{ e\}$.

Now let us consider the second case. By Lemma \ref{b26Mar20}, $n|p^e-1$ (since $\l_n\in \mF_{p^e}$). 
Suppose that  $l:=\gcd( \frac{p^e-1}{n},\gcd_p (g))>1$. Then  $\s_{\l_{ln}, (1-\l_{ln})\nu}\in \oG_f=\langle \s_{\l_n, (1-\l_n)\nu}\rangle $ where $\l_{ln}$ is a primitive $ln$'th root of unity, a contradiction (since the order of the  element $\s_{\l_{ln}, (1-\l_{ln})\nu}$ is $ln>n=|\oG_f|$). Therefore, the statement (a) holds. 

Suppose that the condition (b) does not holds, i.e. there are two distinct roots $\l $ and $\l'$ of the polynomial $f(x)$ such that $v':= \l -\l'\not\in V$ and 
$\CR (f) +\Fp (\l_n)(\l -\l')\subseteq \CR (f)$.  Then $$ V':=V+\Fp (\l_n)v'$$ is a $\Fp (\l_n)$-submodule of $K$ that properly contains the 
$\Fp (\l_n)$-module  $V$. Then $\Sh_{V'} \subseteq \Sh_V$, a contradiction.

$(2\Rightarrow 1)$ In the first case, i.e. $f(x)=f_V^i(x-\nu )$,  the implication follows from Proposition \ref{A26Mar20}. In the second case, i.e.  $f(x)=f_V^i(x-\nu ) g(f_V^n(x-\nu ))$, 
$$\oG_f\supseteq \langle \s_{\l_n, (1-\l_n)\nu}\rangle \neq \{ e\}\;\; {\rm  and }\;\; \widetilde {G_f}\supseteq  \{\s_{1, \mu}\, | \, \mu \in V\}\neq \{ e\}.$$ The conditions (a) and (b) imply that the inclusions above are equalities hold, see the proof of the implication $(1\Rightarrow 2)$. 

Suppose that statement 1 holds. Then statement 2 holds and vice versa. So, we have the equality 
$$ f(x)=f_V^i(x-\nu ) g(f_V^n(x-\nu ))$$ in statement 2 (the case $g=1$ corresponds to the first case). The polynomial $ f_V(x-\nu )$ is $\widetilde{G}_f$-invariant, i.e. for all elements $v\in V$, $f_V(x-\nu)=\s_{1,-v}(f_V(x-\nu ))=f_V(x-(\nu +v))$. Therefore, for all elements $v\in V$, $ f(x)=f_V^i(x-(\nu +v))g(f_V^n(x-(\nu +v)))$, i.e. the element $\nu$ can be replaced by the element $\nu +v$ for any element $v\in V$. This is the only freedom for the choice of the element $\nu$. Indeed, $G_f=\{ \s^j_{\l_n, (1-\l_n)\nu}\s_{1,v}\, | \, 0\leq j\leq n-1, v\in V\}$. Since 
 $$\s^j_{\l_n, (1-\l_n)\nu}\s_{1,v}=\s_{\l_n^j, (1-\l_n^j)\nu}\s_{1,v}=\s_{\l_n^j, (1-\l_n^j)(\nu +(1-\l_n^j)^{-1}v)}\;\; {\rm for}\;\;  1\leq j \leq n-1,$$  it follows that the only freedom in choosing the generator $\s_{\l_n, (1-\l_n)\nu}$ in Theorem \ref{8Mar20} is an element of the type 
   $\s_{\l_n^j, (1-\l_n^j)(\nu +(1-\l_n^j)^{-1}v)}$ where $j$ is a natural number such that $1\leq j \leq n-1$, $\gcd (j,n)=1$ and $v$ is an arbitrary element of $V$. Now, by Theorem \ref{8Mar20}, $\nu$ is unique up to addition an element of $V$, and the elements $i$ and $g(x)$ are unique.
   
    Clearly,  $\s_{\l_n, (1-\l_n)\nu }(f)=\l_n^i f$ (Theorem \ref{8Mar20}.(1)), and so $f\in K[x]^{G_f}$ iff $n|i$  (Theorem \ref{8Mar20}.(2,3)).
   
 Suppose that $g(x) \neq 1$.  Clearly, $\nu$ is a root of $f(x)$ iff $i\geq 1$. Suppose that $\nu$ is not a root of $f(x)$, i.e. $i=0$ and $f(x)=g(f_V^n(x-\nu ))$. Recall that $f(x)=f_1^{p^s}(x)$ and $G_f=G_{f_1}$, by Lemma \ref{a29Mar20}. Then  $\nu$ is not a root of $f_1(x)$, i.e.  $f_1(x)=g_1(f_V^n(x-\nu ))$ for a unique polynomial $g_1(x)\in K[x]$ such that $g=g_1^{p^s}$. Hence, $\nu$ is a root of the polynomial $f_1(x)$ since 
$$0\neq  f_1'(x)=nf_V^{n-1}(x-\nu )f_V'(x-\nu )g_1'(f_V^n(x-\nu )),$$
$n\geq 2$ and $f_V(0)=0$. 

In the second case, i.e.  $f(x)=f_V^i(x-\nu ) g(f_V^n(x-\nu )) $, $n=\gcd (p^e-1, \gcd_p(h)$, by the statement (a).  $\Box$ \\

{\it Definition.} The unique presentation of the polynomial $f(x)$,  
$$f(x)=f_V^i(x-\nu )\;\; {\rm  or}\;\; f(x)=f_V^i(x-\nu ) g(f_V^n(x-\nu )),$$ in Theorem \ref{A10Mar20}.(2) is called the {\em eigenform} or the {\em eigenpresentation} of the polynomial $f(x)$. The scalars $\nu +V$ and the natural number $i\geq 0$ are called the {\em eigenroots} of $f(x)$ and their {\em multiplicity}, respectively. The natural number $n\geq 2$ and the monic polynomial $g(x)$ are called the {\em eigenorder} and the {\em eigenfactor} of $f(x)$. In the second case, the eigenroots may not be roots of the polynomial $f(x)$. They are iff $i\neq 0$.

\begin{corollary}\label{aA10Mar20}
Suppose the field $K$ is an algebraically closed field and    $f(x)\in K[x]$ is a monic  polynomial  that  has at least two distinct roots,  $gcd (f) = p^s\gcd_p(f)$ and $f(x) =f_1^{p^s}(x)$ for a unique monic nonscalar polynomial $f_1(x)\in K[x]$. Suppose that the polynomial $f$ satisfies the assumption of Theorem \ref{A10Mar20}, and $f_1(x)=f_V^j(x-\nu )$ or $ f_1(x)=f_V^j(x-\nu ) g_1(f_V^n(x-\nu ))$, is the eigenform of the polynomial $f_1(x)$. Then $f(x)=f_V^{p^sj}(x-\nu )$ or $ f(x)=f_V^{p^sj}(x-\nu ) g_1^{p^s}(f_V^n(x-\nu ))$, is the eigenform of the polynomial $f_1(x)$.
\end{corollary}

{\it Proof.} The statement follows from the facts that  $f(x)=f_1^{p^s}(x)$,  $G_f=G_{f_1}$ (Lemma \ref{a29Mar20}) and the uniqueness of the eigenform (Theorem \ref{A10Mar20}). $\Box$ \\


{\bf Criterion for $\widetilde{G}_f=\{ e\}$ and $\oG_f\neq \{ e\}$.} Theorem \ref{A11Mar20} is a criterion for the group $G_f=\widetilde{G_f}\rtimes \oG_f$ to be equal to $\oG_f \neq \{ e\}$.

\begin{theorem}\label{A11Mar20}
Suppose the field $K$ is an algebraically closed field,   $f(x)\in K[x]$ is a monic  polynomial  that  has at least two distinct roots,  $gcd (f) = p^s\gcd_p(f)$ and $f(x) =f_1^{p^s}(x)$ for a unique monic nonscalar polynomial $f_1(x)\in K[x]$, see (\ref{fxf1ps}). Then the following statements are equivalent: 
\begin{enumerate}
\item  $\widetilde {G_f}=\{ e\}$ and $\oG_f=\langle \s_{\l_n, (1-\l_n)\nu}\rangle \neq \{ e\}$ where $\l_n$ is a primitive $n$'th root of unity and $\nu \in K$.
\item  $f(x)=(x-\nu )^i g((x-\nu )^n) $ for some natural number $i\geq 0$  and a  monic nonscalar polynomial $g(x)\in K[x]$ such that $g(0)\neq 0$, 
 
\begin{enumerate}
\item $n\geq 2$, $p\nmid n$ and $\gcd_p (g(x))=1$, and 
\item  $\CR (f) +\Fp (\l -\l')\not\subseteq \CR (f)$ for all distinct roots $\l$ and $\l'$ of the polynomial $f$.
\end{enumerate}
\end{enumerate}
Suppose that statement 1 holds. Then the presentation $f(x)=(x-\nu )^i g((x-\nu )^n) $ is unique, i.e. the triple $(\nu , i , g(x))$ is unique. Either $\nu$ is a root of $f(x)$ (i.e. $i\geq 1$) or otherwise (i.e. $i=0$) $\nu$ is a root of $f_1'(x)$ (the derivative of $f_1(x)$).  If  $\nu$ is a root of $f(x)$ then $n=\gcd_p(x^{-i}f(x+\nu ))$.   If $\nu$ is not a root of $f(x)$ then 
 $n=\gcd_p(f(x+\nu ))$. Furthermore, $\s_{\l_n, (1-\l_n)\nu }(f)=\l_n^i f$, and $f\in K[x]^{G_f}$ iff $n|i$.
\end{theorem}

{\it Proof}. By Proposition \ref{A25Mar20}.(3), $\widetilde {G_f}=\{ e\}$ iff the condition (b) holds. 

$(1\Rightarrow 2)$ Suppose that statement 1 holds. Then, by Theorem \ref{A8Mar20}.(2,3), $$f(x)=(x-\nu )^i g((x-\nu )^n)  $$  for some natural number $i\geq 0$ and  a  monic {\em nonscalar} polynomial $g(x)\in K[x]$ such that $g(0)\neq 0$ (since $|\CR_d(f)|\geq 2$). Clearly, $n\geq 2$ and $p\nmid n$. 

Suppose that  $l:=\gcd_p (g(x))>1$. Then  $\s_{\l_{ln}, (1-\l_{ln})\nu}\in \oG_f=\langle \s_{\l_n, (1-\l_n)\nu}\rangle $ where $\l_{ln}$ is a primitive $ln$'th root of unity, a contradiction (since the order of the  element $\s_{\l_{ln}, (1-\l_{ln})\nu}$ is $ln>n=|\oG_f|$). Therefore, the statement (a) holds.

$(2\Rightarrow 1)$  By the statement (b), $\widetilde {G_f}=\{ e\}$. Since  $f(x)=(x-\nu )^i g((x-\nu )^n) $ for some natural number $i\geq 0$, $n\geq 2$, $p\nmid n$   and a  monic nonscalar polynomial $g(x)\in K[x]$ such that $g(0)\neq 0$, 
$$\oG_f\supseteq \langle \s_{\l_n, (1-\l_n)\nu}\rangle \neq \{ e\}.$$ The condition (a)  implies that the inclusion above is the  equality, see the proof of the implication $(1\Rightarrow 2)$.

Suppose that statement 1 holds. Then statement 2 holds and vice versa. So, we have the equality 
$f(x)=(x-\nu )^i g((x-\nu )^n)$ as in statement 2. To prove uniqueness of this presentation it suffices to show that the element $\nu$ is unique.  The set of cyclic generators for the group $G_f=\oG_f=\langle \s_{\l_n, (1-\l_n)\nu}\rangle$ is equal to $\{ \s_{\l_n, (1-\l_n)\nu}^j\, | \, 1\leq j \leq n-1, \gcd (j,n)=1\}$. Since $\s_{\l_n, (1-\l_n)\nu}^j=\s_{\l_n^j, (1-\l_n^j)\nu}$, the element $\nu$ is unique. 

Clearly, $\nu$ is a root of $f(x)$ iff $i\geq 1$. Suppose that $\nu$ is not a root of $f(x)$, i.e. $i=0$ and $f(x)=g((x-\nu )^n)$,  then $f_1(x) =h((x-\nu )^n)$ for a unique monic nonscalar polynomial $h(x):=g^\frac{1}{p^s}(x)\in K[x]$ (since $p\nmid n$). Hence, $\nu$ is a root of the polynomial $f_1(x)$ since 
$$0\neq  f_1'(x)=n(x-\nu )^{n-1}h'((x-\nu )^n)$$
and $n\geq 2$. 

 If  $\nu$ is a root of $f(x)$ then $n=\gcd_p(x^{-i}f(x+\nu ))$ (since $f(x+\nu ) = x^ig(x^n)$).   If $\nu$ is not a root of $f(x)$ then 
 $n=\gcd_p(f(x+\nu ))$ (since $f(x+\nu ) =g(x^n)$). 
 
 Clearly,  $\s_{\l_n, (1-\l_n)\nu }(f)=\l_n^i f$, and so $f\in K[x]^{G_f}$ iff $n|i$. $\Box$ \\

{\it Definition.} The unique presentation of the polynomial $f(x)$,  
$$f(x)=(x-\nu )^i g((x-\nu )^n),$$ in Theorem \ref{A11Mar20}.(2) is called the {\em eigenform} or the {\em eigenpresentation} of the polynomial $f(x)$. The scalar $\nu$ and the natural number $i\geq 0$ are called the {\em eigenroot} of $f(x)$ and its {\em multiplicity}, respectively. The natural number $n\geq 2$ and the monic polynomial $g(x)$ are called the {\em eigenorder} and the {\em eigenfactor} of $f(x)$. In general, the eigenroot may not be a root of the polynomial $f(x)$. It is  iff $i\neq 1$.\\


{\bf Criterion for $\widetilde{G}_f\neq \{ e\}$ and $\oG_f= \{ e\}$.}

\begin{lemma}\label{b28Mar20}
Let $g(x) \in K[x]$ be a monic nonscalar polynomial such that $g'\neq 0$ (the derivative of $g$), and $V$ be an $\Fp$-subspace of $K$ and $\mF_{p^e}$ be its multiplier field (for $V=0$, $f_V(x)=x$). Then 
\begin{enumerate}
\item $g(f_V(x))'\neq 0$.
\item For each $\nu \in K$, $g(f_V(x)) = f_V^{i(\nu )}(x-\nu) g_\nu(f_V(x-\nu ))$ for a unique natural number $i(\nu )\geq 0$ and a unique monic polynomial $g_\nu(x)\in K[x]$ such that $g_\nu(0)\neq 0$;  $i(\nu )\neq 0$ iff $\nu$ is a root of the polynomial $g(f_V(x))$. If 
 $n=\gcd(p^e-1, \gcd_p(g_\nu (x)))\geq 2$ then  $e\neq \s_{\l_n, (1-\l_n)\nu}\in \oG_{g(f_V(x))}(\bK )$ where $\l_n\in \bK$ is a primitive $n$'th root of unity. If, in addition, $\l_n\in K$ then  $e\neq \s_{\l_n, (1-\l_n)\nu}\in \oG_{g(f_V(x))}(K )$.
\item Suppose that $\nu$ is not a root of the polynomial $g(f_V(x))$, i.e.  $i(\nu )= 0$ and  $g(f_V(x))=g_\nu (f_V(x-\nu))$, and $\gcd_p(g_\nu (x))\neq 1$ then $\nu$ is a root of the derivative $g(f_V(x))'$ of the polynomial $g(f_V(x))$. 
\end{enumerate}
\end{lemma}

{\it Proof}. 1. $g(f_V(x))'=g'(f_V(x))f_V'(x)\neq 0$, by Proposition \ref{B7Mar20}.(2d).

2. $g(f_V(x))=g(f_V(x-\nu +\nu))=g(f_V(x-\nu)+f_V(\nu ))=f_V^{i(\nu)}(x-\nu) g_\nu (f_V(x-\nu ))$ for a unique natural number $i(\nu )\geq 0$ and a unique monic polynomial $g_\nu (x) \in K[x]$ such that $ g_\nu (0)\neq 0$. Since $f_V(0)=0$ and $ g_\nu (0)\neq 0$, we see that $i(\nu )\neq 0$ iff $\nu$ is a root of the polynomial $g(f_V(x))$.

If $n\geq 2$ then  $e\neq \s_{\l_n, (1-\l_n)\nu}\in \oG_{g(f_V(x))}(\bK )$, by Theorem \ref{8Mar20}.(1). If, in addition, $\l_n\in K$ then  $e\neq \s_{\l_n, (1-\l_n)\nu}\in \oG_{g(f_V(x))}(K )$.

3. Suppose that $\nu$ is not a root of the polynomial $g(f_V(x))$ and $m=\gcd_p(g_\nu (x))\neq 1$, i.e. $g(f_V(x))=g_\nu (f_V(x-\nu )) =h_\nu (f_V^m(x-\nu ))$ for some monic nonscalar polynomial $h_\nu (x) \in K[x]$. Then 
$$ g(f_V(x))'=h_\nu (f_V^m(x-\nu ))'=mf_V^{m-1}(x-\nu )f'_V(x-\nu ) h_\nu'(f_V^m(x-\nu )),$$ 
and so $\nu$ is a root of the polynomial $g(f_V(x))'$.  $\Box $\\

Given monic nonscalar polynomials $f(x), h(x)\in K[x]$. {\em If $f(x) = g(h(x))$ for some polynomial $g(x)\in K[x]$ then the polynomial $g(x)$ is unique and necessarily monic.} (Proof. If $f(x) = g(h(x))$ then the polynomial $g(x)$ is monic, $\deg (f) = \deg (g)\deg (h)$, $K[h]\ni f_1:= f-h^{\deg (g)}$ and $\deg (f_1)<\deg(f)$. Now, the induction on $\deg (f)$ completes the proof).\\

Theorem \ref{B11Mar20} is a criterion for the group $G_f=\widetilde{G_f}\rtimes \oG_f$ to be equal to $\widetilde{G_f} \neq \{ e\}$.

\begin{theorem}\label{B11Mar20}
Suppose the field $K$ is an algebraically closed field,   $f(x)\in K[x]$ is a monic nonscalar polynomial  that  has at least two distinct roots, $gcd (f) = p^s\gcd_p(f)$ and $f(x) =f_1^{p^s}(x)$ for a unique monic nonscalar polynomial $f_1(x)\in K[x]$. Suppose that  $V\neq 0$ is a  finite dimensional $\Fp $-subspace of $K$ and $\mF_{p^e}$ is the multiplier field of $V$. Then the following statements are equivalent: 
\begin{enumerate}
\item  $\widetilde {G_f}= \Sh_V\neq \{ e\}$ and $\oG_f=\{ e\}$.
\item  $f_1(x) = g(f_V(x))$ for a (unique) monic polynomial $g(x)\in K[x]$ such that 
\begin{enumerate}
\item either $|\CR_d(g)|=1$ and $(p,e)=(2,1)$ or otherwise $|\CR_d(g)|\geq 2$ and $\gcd (p^e-1, \gcd_p(g_{\nu} (x)))=1$ for all roots $\nu \in \CR_d(f_1(x))\cup \CR_d(f_1(x)')$ where $g_\nu$ is as in Lemma \ref{b28Mar20}.(2) (i.e. $f_1(x) = f_V^{i (\nu )}(x-\nu ) g_\nu (f_V(x-\nu ))$),  and 
\item   $\CR (f_1) +\Fp (\l -\l')\not\subseteq \CR (f_1)$ for all distinct roots $\l$ and $\l'$ of the polynomial $f_1$  such that $\l -\l'\not\in V$ ($\Leftrightarrow $  $\CR (f) +\Fp (\l -\l')\not\subseteq \CR (f)$ for all distinct roots $\l$ and $\l'$ of the polynomial $f$  such that $\l -\l'\not\in V$).
\end{enumerate}
\end{enumerate}
Suppose that statement 1 holds. Then $f,f_1\in K[x]^{G_f}=K[x]^{G_{f_1}}$.
\end{theorem}

{\it Remark.} By Lemma \ref{a29Mar20}, $G_f=G_{f_1^{p^s}}=G_{f_1}$. This explains why statement 2 is given via properties of the polynomial $f_1$ rather than $f$.\\

{\it Proof}. By Proposition \ref{B7Mar20} and Proposition \ref{A25Mar20}.(1), 
$\widetilde{G}_{f}=\widetilde{G}_{f_1}=\Sh_V (\neq \{ e\})$ iff  $f_1(x) = g(f_V(x))$ for a (unique) monic nonscalar polynomial $g(x)\in K[x]$ such that the condition 2(b) holds. 

Suppose that $|\CR_d(g)|=1$, i.e. $\CR_d(g)=\{ \rho \}$ and let $i(\rho )$ be the multiplicity of the root $\rho$. Fix an element $\nu'\in K=\bK$ such that $f_V(\nu ')=\rho$. Then $ f_1(x)=(f_V(x)-f_V(\nu '))^{i(\rho )}=f_V^{i(\rho )}(x-\nu')$, by  Proposition \ref{B7Mar20}.(2b).
 By  Proposition \ref{A26Mar20}, $G_{f_1}=\widetilde{G}_{f_1}=\Sh_V$ iff $(p,e)=(2,1)$.

Suppose that $|\CR_d(g)|\geq 2$. By  Lemma  \ref{b28Mar20}.(2), 
 $$ f_1(x) =g(f_V(x))=f_V^{i (\nu )}(x-\nu ) g_\nu (f_V(x-\nu ))$$ for a unique monic  {\em nonscalar} polynomial $g_\nu (x)\in K[x]$ such that $g_\nu (0)\neq 0$ where   $i(\nu )\geq 0$ is the multiplicity of the root $\nu$ (if $g_\nu (x) =1$ then $f_1(x)=g(f_V(x))=f_V^{i (\nu )}(x-\nu )=(f_V(x)-f_V(\nu ))^{i(\nu )}$, and so $|\CR_d(g)|=1$, a contradiction).
 
   By Theorem \ref{8Mar20} and Lemma \ref{b28Mar20}.(2), $\oG_{f_1}=\{ e\}$ iff $\gcd (p^e-1, \gcd_p(g_{\nu} (x)))=1$ for all $\nu \in K$ iff $\gcd (p^e-1, \gcd_p(g_{\nu} (x)))=1$ for all $\nu \in \CR_d(f_1(x))\cup \CR_d(f_1(x)')$, Lemma \ref{b28Mar20}.(2,3). 
   
 Clearly, $f,f_1\in K[x]^{G_f}=K[x]^{G_{f_1}}$ (Proposition \ref{B7Mar20} and Lemma \ref{a29Mar20}).   $\Box$ \\

{\it Definition.} The unique presentation $f(x)=g^{p^s}(f_V(x))$ in Theorem \ref{B11Mar20}
(where $\gcd (f)=p^s\gcd_p(f)$) is called the {\em eigenform} or {\em eigenpresentation} of the polynomial $f(x)$ and the polynomial $g(x)$ is called the {\em eigenfactor} of $f(x)$. \\


{\bf Criterion for $G_f=\{e\}$.}  Given a monic nonscalar polynomial $g(x) \in K[x]$ with $g'(x)\neq 0$. By Lemma \ref{b28Mar20}.(2) (where $V=0$),  for each $\nu \in K$,
\begin{equation}\label{gxxg}
g(x)=(x-\nu )^{i(\nu )}g_\nu (x-\nu )
\end{equation}
for a  natural number $i(\nu )\geq 0$ and a unique monic polynomial $g_\nu (x)\in K[x]$ such that $g_\nu (0)\neq 0$. Clearly, $i(\nu)\neq 0$ iff $\nu$ is a root of the polynomial $g(x)$. \\

Theorem \ref{C11Mar20} is a criterion for $G_f=\{ e\}$.

\begin{theorem}\label{C11Mar20}
Suppose the field $K$ is an algebraically closed field,   $f(x)\in K[x]$ is a monic  polynomial  that  has at least two distinct roots, $\gcd (f)=p^s\gcd_p(f)$ and $f(x)=g^{p^s}(x)$ for a unique monic nonscalar polynomial $g(x)\in K[x]$, see (\ref{fxf1ps}). The following statements are equivalent:
\begin{enumerate}
\item $G_f=\{ e\}$.
\item 
 \begin{enumerate}
  \item  For each root $\nu $ of the polynomial $g(x)$, $\gcd_p(g_\nu (x))=1$ where the polynomial $g_\nu (x)$ is defined in (\ref{gxxg}),
  \item  for each root $\nu' $ of the derivative $g'(x)$ of the polynomial $g(x)$ such that $g(\nu')\neq 0$, $\gcd_p(g_\nu (x))=1$, and 
 
\item   $\CR (g) +\Fp (\l -\l')\not\subseteq \CR (g)$ for all distinct roots $\l$ and $\l'$ of the polynomial $g$ ($\Leftrightarrow $ $\CR (f) +\Fp (\l -\l')\not\subseteq \CR (f)$ for all distinct roots $\l$ and $\l'$ of the polynomial $f$).
\end{enumerate}
\end{enumerate}
\end{theorem}

{\it Proof}. 
 Notice that $G_f=G_{g^{p^s}}=G_g$. The condition (c) is equivalent to the condition that $\widetilde{G}_g=\{ e\}$ (Proposition \ref{A25Mar20}.(3)). It remains to show that provided $\widetilde{G}_g=\{ e\}$ the condition $\oG_g=\{ e\}$ is equivalent to the conditions (a) and (b). Equivalently, $\widetilde{G}_g=\{ e\}$ and  $\oG_g\neq \{ e\}$ iff   one of the conditions (a) or (b) does not hold and the condition (c) holds. This follows  from Theorem \ref{A11Mar20}. Indeed, by Theorem \ref{A11Mar20}, $\widetilde{G}_g=\{ e\}$ and  $\oG_g\neq \{ e\}$ iff the condition (c) holds and $g(x)=(x-\nu )^ig_\nu (x-\nu )$ for a unique $\nu\in K$, a natural number $i\geq 0$ and a monic nonscalar polynomial $g_\nu(x)$ such that $g_\nu (0)\neq 0$ (since $|\CR_d (f)|\geq 2$) and  $\gcd_p(g_\nu (x))\geq 2$. We have two options either $\nu$ is a root of the polynomial $g(x)$ or not. If $\nu$ is not  a root of the polynomial $g(x)$, i.e. $i=0$, then $g(x)=g_\nu (x-\nu )$, and so 
 $\nu$ is a root of $g'(x)$ since  $\gcd_p(g_\nu (x))\geq 2$. Now, it follows that statements 1 and 2 are equivalent.  $\Box$. \\
 

Theorem \ref{A12Mar20} describes the group $G_f(K)$ in terms of the group $G_f(\bK)$.

\begin{theorem}\label{A12Mar20}
Suppose the field $K$ is not necessarily  algebraically closed and    $f(x)\in K[x]$ is a monic  polynomial  that  has at least two distinct roots in $\bK$.  Recall that (Theorem \ref{10Mar20} and Theorem \ref{C7Mar20}) the group $G_f(\bK )=\widetilde{G_f}(\bK )\rtimes \oG_f(\bK )$   where $\oG_f(\bK ):=\langle \s_{\l_n, (1-\l_n)\nu}\rangle $ and $\widetilde {G_f}(\bK ):= \{\s_{1, \mu}\, | \, \mu \in \overline{V}\}$,  $\l_n\in \bK$ is a primitive $n$'th   root of unity provided $\oG_f(\bK )\neq \{ e\}$,  $\nu \in \bK$,  $\overline{V}$ is a  finite dimensional $\Fp (\l_n)$-subspace of $\bK$, and $\Fp (\l_n)=\F_{p^m}$ for some $m\geq 1$ (Lemma \ref{a8Mar20}).  Then $G_f(K )=\widetilde{G_f}(K )\rtimes \oG_f(K )=\Aut_K(K[x])\cap G_f(\bK )$, $\widetilde{G_f}(K )=\Sh_V$ where $V:=K\cap \overline{V}$ and if $\oG_f(K)\neq \{ e\}$ then $\oG_f(K )=\langle \s_{\l_n^i, (1-\l_n^i)\nu +\overline{v}}\rangle $ where $i=\min\{ i'=1, \ldots , n-1\, | \, i'|n, \l_n^{i'}\in K, (1-\l_n^{i'})\nu \in \overline{V}+K\}$ and $\overline{v}\in \overline{V}$ is any (fixed) element such that $(1-\l_n^i)\nu +\overline{v}\in K$. 

\end{theorem}

{\it Proof}.  It is obvious that  $G_f(K )=\Aut_K(K[x])\cap G_f(\bK )$. By Theorem \ref{C7Mar20}, $G_f(K )=\widetilde{G_f}(K )\rtimes \oG_f(K )$.  It is obvious that  $\widetilde{G_f}(K )=\Sh_V$ where $V:=K\cap \overline{V}$. By Theorem \ref{C7Mar20},  $\oG_f(K )=\langle \s_{\l_{n'}, (1-\l_{n'})\nu'}\rangle $ where $\l_{n'}\in K$ is a primitive $n'$'th   root of unity and   $\nu' \in K$ provided $\oG_f(K)\neq \{ e\}$. Notice that
$$\s_{\l_{n'}, (1-\l_{n'})\nu'}=\s^i_{\l_n, (1-\l_n)\nu}\s_{1, \overline{v}}=\s_{\l_n^i, (1-\l_n^i)\nu}\s_{1, \overline{v}}=\s_{\l_n^i, (1-\l_n^i)\nu+\overline{v}}$$ for unique elements $i$ and $\overline{v}\in \overline{V}$ such that $0\leq i\leq n-1$.
So, the elements $i$ can be chosen such that 
\begin{eqnarray*}
 i&=&\min\{ i'=1, \ldots , n-1\, | \, i'|n,\;  \l_n^{i'}\in K, \; (1-\l_n^{i'})\nu+\overline{v}\in K\;\; {\rm  for \; some \;  element }\;\; \overline{v}\in \overline{V}\} \\
 &=&\min\{ i'=1, \ldots , n-1\, | \, i'|n,\;  \l_n^{i'}\in K, \; (1-\l_n^{i'})\nu \in \overline{V}+K\}
\end{eqnarray*}
and $\overline{v}\in \overline{V}$ is any (fixed) element such that $(1-\l_n^i)\nu +\overline{v}\in K$. $\Box$\\

Proposition \ref{a26Mar20} gives criteria for the groups $\widetilde{G}_f(K)$, $\oG_f(K)$ and $G_f(K)$ to be $\{ e\}$. 
 
\begin{proposition}\label{a26Mar20}
Suppose the field $K$ is not necessarily  algebraically closed and    $f(x)\in K[x]$ is a monic  polynomial  that  has at least two distinct roots in $\bK$.  Recall that (Theorem \ref{10Mar20} and Theorem \ref{C7Mar20}) the group $G_f(\bK )=\widetilde{G_f}(\bK )\rtimes \oG_f(\bK )$   where $\oG_f(\bK ):=\langle \s_{\l_n, (1-\l_n)\nu}\rangle $ and $\widetilde {G_f}(\bK ):= \{\s_{1, \mu}\, | \, \mu \in \overline{V}\}$,  $\l_n\in \bK$ is a primitive $n$'th   root of unity provided $\oG_f(\bK )\neq \{ e\}$,  $\nu \in \bK$,  $\overline{V}$ is a  finite dimensional $\Fp (\l_n)$-subspace of $\bK$, and $\Fp (\l_n)=\F_{p^m}$ for some $m\geq 1$ (Lemma \ref{a8Mar20}).  Then 
\begin{enumerate}
\item $\widetilde{G}_f(K)=\{ e\}$ iff $\overline{V}\cap K=0$.
\item $\oG_f(K)=\{ e\}$ iff $\oG_f(\bK)=\{ e\}$ or otherwise $\oG_f(\bK ):=\langle \s_{\l_n, (1-\l_n)\nu}\rangle $ and there is no a natural number $i'$ such that  $1\leq i'\leq n-1$ such that $i'|n$, $\l_n^{i'}\in K$ and $ (1-\l_n^{i'})\nu \in \overline{V}+K$.
\item $G_f(K)=\{ e\}$ iff $G_f(\bK)=\{ e\}$ or otherwise $\overline{V}\cap K=0$, $\oG_f(\bK ):=\langle \s_{\l_n, (1-\l_n)\nu}\rangle $ and there is no a natural number $i'$ such that  $1\leq i'\leq n-1$ such that $i'|n$, $\l_n^{i'}\in K$ and $ (1-\l_n^{i'})\nu \in \overline{V}+K$.
\end{enumerate}
\end{proposition}

{\it Proof}. Statements 1 and 2 follow at once from Theorem \ref{A12Mar20}. Then statement 3 follows from statements 1 and 2. $\Box $\\

{\bf Every subgroup of $\Aut_K(K[x])$ is of type $G_f$.}  Theorem \ref{30Mar20} shows that all subgroups of  $\Aut_K(K[x])$ are eigengroups of polynomials.
 
\begin{theorem}\label{30Mar20}
Let $K$ be an arbitrary field of characteristic $p>0$. Then for each subgroup $H$ of  $\Aut_K(K[x])$ there is a monic polynomial $f_H$ such that $G_{f_H}=H$: 
\begin{enumerate}
\item For $H=\{ e\}$, $f_H=x(x+1)^2$.
\item For $H=\langle \s_{\l_n, (1-\l_n)\nu}\rangle$ where $\l_n\in K$ is a primitive $n$'th root of unity and $\nu \in K$, $f_H=(x-\nu)^n-1$. 
\item For $H=\Sh_V$ where $V$ is a nonzero $\Fp$-subspace of $K$, 
\begin{enumerate}
\item if $K=\mF_{p^n}$ then $f_H(x)=f_V(x-\nu )-\rho$ where $\rho$ is any element of $\mF_{p^n}$ that does not belong to the image of the map $f_V(x-\nu ):\mF_{p^n}\ra \mF_{p^n}$, $x\mapsto f_V(x-\nu )$ (the map $f_V(x-\nu )$ is not a surjection since the set $\nu +V$ is mapped to $0$). 
\item If $|K|=\infty$ then $f_H(x) = f_V(x) f_V^2(x-\nu )$ where $\nu \in K\backslash V$. 
\end{enumerate}
\item For $H=\Sh_V\rtimes \langle \s_{\l_n, (1-\l_n)\nu}\rangle$ where $V$ is a nonzero $\Fp$-subspace of $K$, $\mF_{p^e}$ is its multiplier field, and $\l_n$ is primitive $n$'th root of unity such that $\l_nV\subseteq V$, $f_H(x)=\begin{cases}
f_V(x-\nu )& \text{if }n=p^e-1,\\
f_V^n(x-\nu )+1& \text{if }n<p^e-1.\\
\end{cases}
$
\item For $H=\mT_\nu (K)=\{  \s_{\l , (1-\l )\nu}
\, | \, \l \in K^\times \}$, $f_H =x-\nu$.
\item For $H=\Aut_K(K[x])$, $f_H=\begin{cases}
1& \text{if }|K|=\infty,\\
x^{p^n}-x& \text{if }K=\mF_{p^n}.\\
\end{cases}
$
\end{enumerate}
\end{theorem}

{\it Proof}. 1. If $\s \in G_{f_H}$ then the maximal ideals $(x)$ and $(x+1)$ of $K[x]$ are $\s$-stable, hence $\s =e$.

2. Clearly, $G_{f_H}\supseteq H$. 

(i) $\widetilde{G}_{f_H}=\{ e\}$: The polynomial $f_H$ has $n$ distinct roots, namely, $\{\nu+\l_n^i\, | \, i=0,1,\ldots , n-1\}$, and $p\nmid n$. Suppose that $\widetilde{G}_{f_H}=\Sh_V\neq \{ e\}$ for some nonzero $\Fp$-subspace $V$ of $K$. Then $p||V|$. Since $V+\CR (f)\subseteq \CR (f)$, we must have $|V|||\CR (f)|$, i.e. $|V||n$, and so $p|n$, a contradiction.

(ii) $G_{f_H}=H$: Let $\s =\s_{\l_n, (1-\l_n)\nu}$. By the statement (i),  $G_{f_H}=\oG_{f_H}=\langle \s'\rangle$ where $\s =\s_{\l_m, (1-\l_m)\nu'}$ for some primitive $m$'th root of unity $\l_m$  and $\nu'\in K$. Since $H\subseteq \oG_{f_H}$, we must have $\nu'=\nu$ and $ n|m$ (since $(x-\nu')$ is the only $\langle \s'\rangle$-invariant maximal ideal of $K[x]$, $(x-\nu)$ is the only $\langle \s\rangle$-invariant maximal ideal of $K[x]$ and $\langle \s\rangle\subseteq \langle \s'\rangle$). 

The polynomial $f_H$ is an eigenvector for the automorphism $\s'$ with (necessarily) eigenvalue $\l_m^n$ since $ n=\deg (f_H)$. Now, $\l_m^nf_H=\l_m^n((x-\nu )^n-1) =\s' (f_H)= \l_m^n (x-\nu )^n-1$. Therefore, $\l_m^n=1$, and so 
$\langle \s\rangle = \langle \s'\rangle$. This means that $\oG_{f_H}=H$, and the statement (ii) follows from the statement (i).

3(a). Clearly, $G_{f_H}\supseteq H =\Sh_V$. 

(i) $\widetilde{G}_{f_H}=H$: The statement follows at once from the fact that the polynomial $f_H$ has $|V|$ distinct roots in $\bK$ (if $\widetilde{G}_{f_H}=\Sh_{V'}$ for some $\Fp$-subspace $V'$ of $K$ that properly contains $V$ then the polynomial $f_H$ contains at least $|V'|$ distinct roots in $\bK$, a contradiction). 

(ii) $\oG_{f_H}=\{ e\}$: Suppose that $\oG_{f_H}\neq \{ e\}$. Then $\oG_{f_H}=\langle \s \rangle$ where $\s = \s_{\l_n, (1-\l_n)\nu'}$, $1\neq \l_n\in K$ is a primitive $n$'th root of unity such that $\l_nV\subseteq V$  and $\nu'\in K$. Notice that $\s (f_H)=\l_n^{|V|}f_H$ and 
$$ f_H(x) = f_V(x-\nu'-(\nu -\nu'))-\rho= f_V(x-\nu')-f_V(\nu -\nu')-\rho=f_V(x-\nu')+f_V(\nu'-\nu )-\rho .$$
By Theorem \ref{8Mar20}.(1), $\s (f_V(x-\nu'))=\l_nf_V(x-\nu')$. Let $a=f_V(\nu'-\nu )-\rho$. Now,
$$\l_n^{|V|}(f_V(x-\nu') +a)=\l_n^{|V|}f_H=\s (f_H)=\s (f_V(x-\nu') +a)=\l_nf_V(x-\nu')+a.$$
Hence, $\l_n^{|V|}=\l_n\neq 1$ and $(\l_n-1)a=0$, i.e. $a=0$. The last equality implies that $\rho\in \im \, f_V(x-\nu)$, a contradiction, and the statement (ii) follows. 

(b). Clearly, $G_{f_H}\supseteq H=\Sh_V$. 

(i)  $\widetilde{G}_{f_H}=H$: The statement follows at once from the fact that $\CR (f) =V\coprod (\nu+V)^2$ where the upper index $`2'$ means that the multiplicity of each root in $\nu +V$ is $2$. 

(ii) $\oG_{f_H}=\{ e\}$: Suppose that $e\neq \s \in \oG_{f_H}$. Then $\s \in \oG_{f_V(x)}\cap \oG_{f_V(x-\nu )}$. Notice that $\oG_{f_V(x-\nu )}=\langle \s_{\l_n, (1-\l_n)\nu}\rangle$ for some primitive $n$'th root of unity $\l_n$ such that $\Fp (\l_n)V\subseteq V$. Then there is a natural number $i$ such that $\s = \s^i_{\l_n, (1-\l_n)\nu}=\s_{\l_n^i, (1-\l_n^i)\nu}\in \oG_{f_V(x)}$. In particular, $V\ni \s^{-1} *(0)=(1-\l_n^i)\nu $, and so $\nu \in (1-\l_n^i)^{-1}V=V$, a contradiction ($1-\l_n^i\neq 0$ since $\s \neq e$). 

4. Statement 4 follows from Theorem \ref{A10Mar20}.

5 and 6. Statements 5 and 6 are obvious. 
$\Box $\\

{\bf Algorithm of finding the eigengroup $G_f(\bK)$ and the eigenform of $f$.} 
The algorithm consists of finitely many steps and is based on Prpoposition \ref{A9Mar20}, Theorem \ref{A10Mar20}, Theorem \ref{A11Mar20}, Theorem \ref{B11Mar20} and Theorem \ref{C11Mar20}.
We assume that $K=\bK$.\\

{\em Step 1.} If $|\CR_d(f)|=1$ then apply Proposition \ref{A9Mar20} to find $G_f$.

From this moment on we assume that  $|\CR_d(f)|\geq 2$.

{\em Step 2.} Use Theorem \ref{C11Mar20} to check whether $G_f=\{e\}$ or $G_f\neq \{e\}$.

From this moment on we assume that $G_f\neq \{e\}$.

{\em Step 3.} By Proposition  \ref{A25Mar20}.(1), the group $\widetilde{G}_f=\Sh_V$ can be found. 

{\em Step 4.} Suppose that $\widetilde{G}_f=\{e\}$. Then necessarily $\oG_f\neq \{ e\}$, and using Theorem \ref{A11Mar20} the group $\oG_f$ is found. In more detail, we know that $\oG_f=\langle \s_{\l_n, (1-\l_n)\nu}\rangle $ and  that $f(x)=(x-\nu )^ig((x-\nu )^n)$ for a unique $\nu \in \CR_d(f)\cup  \CR_d(f_1')$ and $n\geq 2$ such that if  $\nu$ is a root of $f(x)$ then $n=\gcd_p(x^{-i}f(x+\nu ))$, and    if $\nu$ is not a root of $f(x)$ then  $n=\gcd_p(f(x+\nu ))$.

From this moment on we assume that $\widetilde{G}_f=\Sh_V\neq \{ e\}$, the $\Fp$-subspace $V$ of the field $K$ is non-zero.
 Let $\mF_{p^e}$ be the multiplier field of $V$. It can be easily found since the multiplier  field $\mF_{p^e}$ is the largest among  finite fields $\mF_{p^m}$ such  that $\mF_{p^m}V\subseteq V$ and $m\leq \dim_{\Fp} (V)$.

{\em Step 5.} Now, we check whether the conditions of Theorem \ref{B11Mar20} hold or not. If they do then $\oG_f=\{ e\}$.

If they do not then necessarily $\oG_f \neq \{e\}$ and hence the conditions of Theorem  Theorem \ref{A10Mar20} hold. Using  Theorem \ref{A10Mar20} the group $\oG_f$ and the eigenform of $f$ are found in finitely many steps.  $\Box$\\

{\bf Algorithm of finding the eigengroup $G_f(K)$ where $K\neq \bK$.}

{\em Step 1.} Using the algorithm above 
 the group $G_f(\bK )$ is found.

{\em Step 2.} The group $G_f(K )$ is found by using Theorem \ref{A12Mar20}. $\Box$

$${\bf Acnowledgements} $$

The author would like to thank the Royal Society for support.

Department of Pure Mathematics

University of Sheffield

Hicks Building

Sheffield S3 7RH

UK

email: v.bavula@sheffield.ac.uk

\end{document}